\documentclass[a4paper]{article}

\usepackage[english]{babel}
\usepackage[utf8x]{inputenc}
\usepackage[T1]{fontenc}
\usepackage{amsmath,amsthm}
\usepackage{amsmath}
\usepackage{amsfonts}
\usepackage{amssymb}
\usepackage{amsfonts,amsmath,amsthm,amscd,amssymb,latexsym}

\usepackage[a4paper,top=3cm,bottom=2cm,left=3cm,right=3cm,marginparwidth=1.75cm]{geometry}

\usepackage{stackrel}
\usepackage{xargs}

\usepackage{amsmath}
\usepackage{graphicx}
\usepackage[colorinlistoftodos]{todonotes}
\usepackage[colorlinks=true, allcolors=blue]{hyperref}

\newtheorem{theorem}{Theorem}[section]
\newtheorem{lemma}[theorem]{Lemma}
\newtheorem{proposition}[theorem]{Proposition}
\newtheorem{corollary}[theorem]{Corollary}

\newtheorem{definition}[theorem]{Definition}

\newtheorem{remark}[theorem]{Remark}
\providecommand{\keywords}[1]{\textbf{\textit{Key words and phrases: }} #1}
\providecommand{\MSC}[1]{\textbf{\textit{2010 MSC: }} #1}

\def\Z{\mathbb{Z}}

\def\N{\mathbb{N}}
\def\K{\mathbb{K}}

\def\A{A}
\def\O{\mathcal{O}_{q,\delta}^{r,n}}
\def\o{\mathcal{O}_{q,\delta}}
\def\u{\bar{u}}
\def\U{\tilde{u}}

\def\pa{g_s}
\def\Pa{g_s}

\newcommandx\p[1]{\left(  #1 \right) }
\newcommandx\br[1]{\lbrace  #1 \rbrace }

\title{Poisson brackets on some skew PBW extensions.}
\author{  Brian Andres Zambrano Luna\footnote{e-mail: bazambranol@unal.edu.co.}\\
Seminario de Álgebra Constructiva - $SAC^2$\\
Departamento de Matemáticas, Universidad Nacional de Colombia, Sede Bogotá.\\ \textit{This article is dedicated to my family.}}
\date{May 5, 2020}

\begin{document}
\maketitle

\begin{abstract}
In [1] the author gives a description of Poisson brackets on  some algebras of quantum polynomials $\mathcal{O}_q$, which is called\textit{ general algebra of quantum polynomials}. The main of this paper is to present a generalization of [1] through a description of Poisson brackets on some skew PBW extensions of a ring $A$ by the extensions $\O$, which are generalization of $\mathcal{O}_q$, and show some examples of skew PBW extension where we can apply this description.   
\end{abstract}
\MSC{16S10 }

\section{Introduction.}
Skew PBW extensions were introduced in [2] and some of its properties have been studied in [4], [5], [6], [7], among others. Let A and R be rings, we say that R is a skew PBW extension of A, if the following conditions hold:
\begin{enumerate}
\item $A\subset R$
\item There exist finitely many elements $x_1,\cdots,x_n\in R$ such that $R$ is a left A-free module with basis $Mon\{x_1,\cdots,x_n\}=\{x_1^{\alpha_1}\cdots x_n^{\alpha_n}|(\alpha_1,\cdots,\alpha_n)\in\N\}$.
\item For every $i=1,\cdots,n$ there exist an injective ring endomorphism $\sigma_i:A\rightarrow A$ and a $\sigma_i$-derivation $\delta_i:A\rightarrow A$ such that 
\[
x_ia=\sigma_i(a)x_i+\delta_i(a)
\]
for all $a\in A$.
\item For every $1\leq i<j\leq n $ there exists $q_{ij}\in A$ left invertible and $a_{ij}^{(t)}\in A$ such that 
\[
x_jx_i=q_{ij}x_ix_j+a_{ij}^{(0)}+a_{ij}^{(1)}x_1+\cdots+a_{ij}^{(n)}x_n
\]
\end{enumerate}
Under these conditions we will denote $R=\sigma(A)\langle x_1,\cdots,x_n\rangle$.

Fix $1 \leq r \leq n$, we will denote  $\O$ an extension of A such that

\begin{enumerate}
\item $A\subset \O$.
\item There exist finitely many elements $x_1,\cdots,x_n,x_1^{-1},\cdots,x_r^{-1}\in \O$ such that $\O$ is a left A-free module with basis $Mon\{x_1^{\pm},\cdots,x_r^{\pm},x_{r+1},\cdots,x_n\}=\{x_1^{\alpha_1}\cdots x_n^{\alpha_n}\vert (\alpha_1,\cdots,\alpha)\in Z\}$, where $Z=\Z^r\times \N^{n-r}$.  
\item For every $i=1,\cdots,n$ there exist a derivation $\delta_i:A\rightarrow A$ such that 
\[
x_ia=ax_i+\delta_i(a)
\]
for all $a\in A$.
\item For every $1\leq i<j\leq n $ there exists $q_{ij}\in Z(A)$  invertible and $a_{ij}^{(t)}\in Z(A)$ such that 
\[
x_jx_i=q_{ij}x_ix_j+a_{ij}^{(0)}+a_{ij}^{(1)}x_1+\cdots+a_{ij}^{(n)}x_n
\]
\item $x_ix_i^{-1}=x_i^{-1}x_i=1$ for $i=1,\cdots,r$.
\end{enumerate}
In the present paper we shall assume that $\O$ satisfies the following conditions  
\begin{enumerate}
\item $\delta_1= 0$
\item For every $i$ fixed, with $i=1,\cdots,n$ and $(m_1,\cdots,m_n)\in Z\setminus \{(0,\cdots,0)\}$, $\p{1-\prod_{j=1,j\neq i}^n q_{ij}^{m_j} }\in \A^{\ast}$.
\item $\delta_t(q_{ij})=\delta_t(a_{ij}^{(m)}) =0$ for all $m = 0, \cdots, n$ and  $i,j,t = 1, \cdots, n$.
\item $a_{1j}^{(0)}=0 $ for all $j=1,\cdots,n$. If $1\leq i\leq r $ and $i<j$ then   $a_{ij}^{(m)}=0$ for all $m=i,\cdots ,n$.
\end{enumerate}
We will denote $\bar{p}_{ij}=p_{ij}-a_{ij}^{(0)}$.
\remark We will consider the deglex order over $Mon(\O):=Mon\{x_1^{\pm},\cdots,x_r^{\pm},x_{r+1},\cdots,x_n\}$ that is defined by [2] as
\[x^{\alpha}\succeq x^{\beta}= \left \{ \begin{array}{lll}
x^{\alpha}= x^{\beta}   &   \\
           &      or         &\\ 
x^{\alpha}\neq x^{\beta} \, but\, \mid \alpha \mid > \mid \beta \mid   &     \\
    &        or       &\\ 
  x^{\alpha}\neq x^{\beta}, \mid \alpha \mid =\mid \beta \mid \, but\, \,\exists i\, \, with\, \alpha_1=\beta_i,\cdots,\alpha_{i-1}=\beta_{i-1} \, and\, \alpha_i>\beta_i   &   \\
         
     \end{array}
\right. \]
where $x^{\alpha}, x^{\beta} \in Mon(\O)$. Each element $f\in \O\setminus\{0\}$ can be represented in a unique way as $f=\eta_{v_{0}}x^{v_0}+\cdots+\eta_{v_k}x^{v_k}$, with $\eta_{v_l}\in A\setminus\{0\}$, $1\leq l\leq k$, and $x^{v_k}\succ\cdots\succ x^{v_0}$; we take $deg(f)=deg(x^{v_k}):=\mid v_k \mid$\footnote{If $\alpha=(\alpha_1,\cdots,\alpha_n)\in\Z^n$ then we take $\mid\alpha\mid=\alpha_1+\cdots+\alpha_n$}. We will denote $it(f)=\eta_{v_k}  x^{v_k}$ the \it{leader term} of f. 

\lemma For all $z\in\O$, $i,j=1,\cdots,n$, and $m=0,\cdots,n$ we have that $q_{ij}z=zq_{ij}$ and $a_{ij}^{(m)}z=za_{ij}^{(m)}$.

\begin{proof}
It is enough to see that $q_{ij}x_t=x_tq_{ij}$ and $a_{ij}^{(m)}x_t=x_ta_{ij}^{(m)}$  for all $t$. By the definition of $\O$ we have that $x_tq_{ij}=q_{ij}x_t+\delta_t(q_{ij})=q_{ij}x_t$ and $x_ta_{ij}^{(m)}=a_{ij}^{(m)}x_t+\delta_t(a_{ij}^{(m)})=a_{ij}^{(m)}x_t$. 
\end{proof}

\lemma Let $x_1^{a_1}\cdots x_n^{a_n}\in Mon(\O)$ then 
\begin{enumerate}
\item For all $i=1,\cdots,n$ there exists $p_{i,a}\in\O$ such that $x_1^{a_1}\cdots x_n^{a_n}x_i=\p{\prod_{j>i}q_{ij}^{a_j}}x_1^{a_1}\cdots x_i^{a_i+1}\cdots x_n^{a_n} +p_{i,a}$ with $deg(p_{i,a})<a_1+\cdots +a_n+1$ or $p_{i,a}=0$ where $a=(a_1,\cdots,a_n)$.
\item For all $i=1,\cdots,n$ there exists $p_{a,i}\in\O$ such that $x_ix_1^{a_1}\cdots x_n^{a_n}=\p{\prod_{j<i}q_{ji}^{a_j}}x_1^{a_1}\cdots x_i^{a_i+1}\cdots x_n^{a_n} +p_{a,i}$ with $deg(p_{a,i})<a_1+\cdots+a_n+1 $ or $p_{a,i}=0$ where $a=(a_1,\cdots,a_n)$. 
\end{enumerate}

\begin{proof}
Fix $x_i\in \O$ and take $x_{t_1}^{a_1}\cdots x_{t_j}^{a_j}\in\O$ such that $a_l\neq0$. We will prove the first claimed by induction on $j$, 
\begin{enumerate}
\item ($j=1, x_{t_1}^{a_1}=x_t^b$)   We will show this claimed by induction on $b$, we can suppose that $i<t$, in other way, we have the claimed.
\begin{enumerate}
\item ($b=1$) By definition, for all $i<t$ we have that $x_tx_i=q_{it}x_ix_t+p_{i,t}$, where $deg(p_{i,t})<2$. 
\item ($b+1$) Put $i<t $. By induction hypothesis there exist $p_{l,b} $ such  that $x_t^bx_l=q_{lt}^bx_lx_t^b+p_{l,b}$ where $deg(p_{l,b})<b+1$ for all $l<t$, then
\begin{eqnarray*}
x_t^{b+1}x_i &=& x_t^b\p{q_{it}x_ix_t+a_{it}^{(0)}+\sum_{l=1}^na_{it}^{(l)}x_l}\\[0.15cm]
&=&q_{it}x_t^{b}x_ix_t+\p{a_{it}^{(0)}x_t^b+\sum_{l=1}^na_{it}^{(l)}x_t^bx_l} \\[0.15cm]
&=& q_{it}q_{it}^bx_ix_t^{b+1}+q_{it}p_{i,b}+a_{it}^{(0)}x_t^b+\sum_{l=1}^{t-1}a_{it}^{(l)}q_{lt}^bx_lx_t^b+\sum_{l=1}^{t-1}a_{it}^{(l)}p_{l,b}+\sum_{l=t}^na_{it}^{(l)}x_t^bx_l \\[0.15cm]
&=&q_{it}^{b+1}x_ix_t^{b+1}+p_{i,b+1}
\end{eqnarray*}
where $ p_{i,b+1}=a_{it}^{(0)}x_t^b+\sum_{l=1}^{t-1}a_{it}^{(l)}q_{lt}^bx_lx_t^b+\sum_{l=1}^{t-1}a_{it}^{(l)}p_{l,b}+\sum_{l=t}^n a_{it}^{(l)}x_t^bx_l $ and

\begin{equation*}
deg(a_{it}^{(0)}x_t^b ),\; deg(a_{it}^{(l)}q_{lt}^bx_lx_t^b ),\; deg(a_{it}^{(l)}p_{l,b} ),\; deg( a_{it}^{(l)}x_t^bx_l)<b+2
\end{equation*}
\item ($b=-1$, if $t\leq r$) We will show that $p_{i,-1}=\sum_{h=1}^{i-1}\sum_{k=1}^{s(i,-1)} r_{h,k}x_{h}x_t^{-k}$ with $r_{h,k}\in\langle a_{ij}^{(m)},q_{ij}| i,j=1,\cdots, n,\, m=0,\cdots,n   \rangle=G$,  sub-ring of $A$, by induction on $1\leq i\leq r$. If $i=1$ then $x_t^{-1}x_1=x_t^{-1}\p{x_1x_t}x_t^{-1}=x_t^{-1}\p{q_{t1}x_tx_1}x_t^{-1}=q_{1t}^{-1}x_1x_t^{-1} $ and $p_{1,-t}=0$. If $i>1$ by induction hypothesis we have that  
\begin{eqnarray*}
x_{t}^{-1}x_i&=&x_t^{-1}\p{x_ix_t}x_t^{-1} =x_{t}^{-1}\p{q_{ti}x_{t}x_i+ a_{ti}^{(0)}+\sum_{l=1}^{i-1}a_{ti}^{(l)}x_l }x_{t}^{-1}\\[0.15cm]
&=& q_{ti}x_ix_{t}^{-1}+a_{ti}^{(0)}x_{t}^{-2}+\sum_{l=1}^{i-1}a_{ti}^{(l)}x_{t}^{-1}x_lx_{t}^{-1} \\[0.15cm]
&=& q_{it}^{-1}x_ix_{t}^{-1}+a_{ti}^{(0)}x_{t}^{-2}+\sum_{l=1}^{i-1}a_{ti}^{(l)}q_{lt}^{-1}x_lx_{t}^{-2}+\sum_{l=1}^{i-1}a_{ti}^{(l)} p_{l,-t}x_t^{-1}
\end{eqnarray*}
since $p_{l,-t}=\sum_{h=1}^{l-1}\sum_{k=1}^{s(l,-1)} r_{h,k}x_{h}x_t^{-k}$ with  $deg(p_{l,-t})<0$ for all $l=1,\cdots,i-1$ then $deg(p_{l,-t})x_t^{-1}<-1$ for all $l$ and $p_{l,-t}x_t^{-1}=\sum_{h=1}^{l-1}\sum_{k=1}^{s(l,-1)} r_{h,k}x_{h}x_t^{-k-1} $ so we have the claimed.
\item ($-(b+1)$, $b>0$ if $t\leq r$) The same way, put $i<t$,  by induction hypothesis there exist $p_{i,-1},p_{i,-b}$ as in  (c) such that $x_t^{-1}x_i=q_{it}^{-1}x_ix_t^{-1}+p_{i,-1} $ and $x_t^{-b}x_i=q_{i,-t}^{-b}x_ix_t^{-b}+p_{i,-b} $, with $deg(p_{i,-1})<0  $ and $deg(p_{i,-b})<-b+1$ then 
\begin{eqnarray*}
x_t^{-b-1}x_i &=& x_t^{-1}(x_t^{-b}x_i) \\[0.15cm]
&=& x_t^{-1}\p{q_{it}^{-b}x_ix_t^{-b}+p_{i,-b}} \\[0.15cm]
&=&q_{it}^{-b}x_t^{-1}x_ix_t^{-b}+x_t^{-1}p_{i,-b}\\[0.15cm]
&=& q_{it}^{-b}\p{q_{it}^{-1}x_ix_t^{-1}+b_{i,-1}}x_t^{-b}+x_t^{-1}p_{i,-b} \\[0.15cm]
&=&q_{it}^{-b-1}x_ix_t^{-b-1}+q_{it}^{-b}b_{i,-1}x_t^{-b}+x_t^{-1}\sum_{h=1}^{i-1}\sum_{k=1}^{s(i,-b)} r_{h,k}x_{h}x_t^{-k}\\[0.15cm]
&=& q_{it}^{-b-1}x_ix_t^{-b-1}+q_{it}^{-b}b_{i,-1}x_t^{-b}+\sum_{h=1}^{i-1}\sum_{k=1}^{s(i,-b)} r_{h,k}x_t^{-1}x_{h}x_t^{-k}\\[0.15cm]
&=&q_{it}^{-b-1}x_ix_t^{-b-1}+q_{it}^{-b}b_{i,-1}x_t^{-b} \\[0.15cm]
&+&\sum_{h=1}^{i-1}\sum_{k=1}^{s(i,-b)} r_{h,k}q_{ht}^{-1}x_hx_t^{-b-1}+\sum_{h=1}^{i-1}\sum_{k=1}^{s(i,-b)} r_{h,k}p_{h,-1}x_t^{-b}
\end{eqnarray*}
Note that $x_tr_{h,k}=r_{h,k}x_t $ for all $h,k$ because for all $t$ and $g\in G$, $\delta_t(g)=0$  so
\[
p_{i,-b-1}=b_{i,-1}x_t^{-b}+\sum_{h=1}^{i-1}\sum_{k=1}^{s(i,-b)} r_{h,k}q_{ht}^{-1}x_hx_t^{-b-1}+\sum_{h=1}^{i-1}\sum_{k=1}^{s(i,-b)} r_{h,k}p_{h,-1}x_t^{-b}
\]
where we have the claimed.
\end{enumerate}
\item $(j+1)$ We can suppose that $t_i=i$ and $i<t_{j+1}$  because in other way we can put $x_i^{a_i}$ with $a_i=0$.  By induction hypothesis on $j=1$ there exist $p_{i,t_{j+1}}$ such that $deg(p_{i,t_{j+1}})<a_{j+1}+1$ and $x_{t_{j+1}}^{a_{j+1}}x_i=q_{it_{j+1}}^{a_{j+1}}x_i x_{t_{j+1}}^{a_{j+1}} + p_{i,t_{j+1}}$, again by induction hypothesis  there exist $p_{i,t'}$, $t'=(t_1,\cdots,t_j) $, such that $ deg(p_{i,t'})< a_1+\cdots+a_i+1$ and $x_{t_1}^{a_1}\cdots x_{t_i}^{a_i}\cdots x_{t_j}^{a_j}x_i =\p{\prod_{l=i+1}^{j}q_{it_l}^{a_l}}x_{t_1}^{a_1}\cdots x_{t_i}^{a_i+1}\cdots x_{t_j}^{a_j}+p_{i,t'} $ then 
\begin{eqnarray*}
x_{t_1}^{a_1}\cdots x_{t_j}^{a_j}x_{t_{j+1}}^{a_{j+1}}x_i &=&x_{t_1}^{a_1}\cdots x_{t_j}^{a_j}\p{q_{it_{j+1}}^{a_{j+1}}x_i x_{t_{j+1}}^{a_{j+1}} + p_{i,t_{j+1}}} \\[0.15cm]
&=& \p{x_{t_1}^{a_1}\cdots x_{t_j}^{a_j}q_{it_{j+1}}^{a_{j+1}}x_i x_{t_{j+1}}^{a_{j+1}} }+\p{x_{t_1}^{a_1}\cdots x_{t_j}^{a_j}p_{i,t_{j+1}} } \\[0.15cm]
&=&q_{it_{j+1}}^{a_{j+1}}\p{x_{t_1}^{a_1}\cdots x_{t_j}^{a_j}x_i x_{t_{j+1}}^{a_{j+1}} }+\p{x_{t_1}^{a_1}\cdots x_{t_j}^{a_j}p_{i,t_{j+1}} }\\[0.15cm]
&=& q_{it_{j+1}}^{a_{j+1}} \p{\p{\prod_{l=i+1}^{j}q_{it_l}^{a_l}}x_{t_1}^{a_1}\cdots x_{t_i}^{a_i+1}\cdots x_{t_j}^{a_j}+p_{i,t'} }x_{t_{j+1}}^{a_{j+1}} \\[0.15cm]
&+&\p{x_{t_1}^{a_1}\cdots x_{t_j}^{a_j}p_{i,t_{j+1}} }\\[0.15cm]
&=& \p{\prod_{l=i+1}^{j+1}q_{it_l}^{a_l}}x_{t_1}^{a_1}\cdots x_{t_i}^{a_i+1}\cdots x_{t_{j+1}}^{a_{j+1}} \\[0.15cm]
&+& q_{it_{j+1}}^{a_{j+1}}p_{i,t'}x_{t_{j+1}}^{a_{j+1}}+ x_{t_1}^{a_1}\cdots x_{t_j}^{a_j}p_{i,t_{j+1}}\\[0.15cm]
&=& \p{\prod_{l=i+1}^{j+1}q_{it_l}^{a_l}}x_{t_1}^{a_1}\cdots x_{t_i}^{a_i+1}\cdots x_{t_{j+1}}^{a_{j+1}} + p_{i,t}
\end{eqnarray*}
where $deg(q_{it_{j+1}}^{a_{j+1}}p_{i,t'}x_{t_{j+1}}^{a_{j+1}} )\leq deg(p_{i,t'} )+deg( x_{t_{j+1}}^{a_{j+1}})<a_1+\cdots+a_j+1+a_{j+1}$ and $deg(x_{t_1}^{a_1}\cdots x_{t_j}^{a_j}p_{i,t_{j+1}} )\leq deg(x_{t_1}^{a_1}\cdots x_{t_j}^{a_j} )+deg(p_{i,t_{j+1}} )<a_1+\cdots+a_j+a_{j+1}+1$ later $deg(p_{it})=deg(q_{it_{j+1}}^{a_{j+1}}p_{i,t'}x_{t_{j+1}}^{a_{j+1}}+ x_{t_1}^{a_1}\cdots x_{t_j}^{a_j}p_{i,t_{j+1}})<a_1+\cdots+a_{j+1}+1$. The second claimed is proved the same way we proved the first claimed.
\end{enumerate}
\end{proof}

\section{Derivations over $\O$.}
\definition Let $d$ be a linear operator over $\O$ such that 
\begin{equation*}
d(ab)=d(a)b+\gamma(a)d(b)+\sum_{s}\theta_s\alpha_s(a)\beta_s(b)
\end{equation*}
for all $a,b\in \O$, where $\theta_s\in\A$ and $\gamma,\alpha_s,\beta_s$ are toric automorphisms of $\A$, i.e. $\gamma(x_i)=\gamma_ix_i$, $\alpha_s(x_i)=a_{s,i}x_i$, and $\beta_s(x_i)=b_{s,i}x_i$ for all $i=1,\cdots,n$ where $\gamma_i,a_{s,i},b_{s,i}\in\A$.

\begin{enumerate}
\item If $\theta_s=0$ for all $s$ then $d$ is a $\gamma$-derivation.
\item If $\theta_s=0$ and $\gamma$ is the identity of $\A$ then $d$ is a derivation.
\item A inner $\gamma$-derivation $[ad_{\gamma}u]a$ is defined by 
\begin{equation*}
[ad_{\gamma}u]a:=ua-\gamma(a)u
\end{equation*}
when $\gamma$ is the identity we denote $[ad_{\gamma}u]a=[ad\, u]a$.
\end{enumerate}

\lemma Let $u_i=d(x_i)$ then they are solution of 

\begin{equation}
u_jx_i+\gamma(x_j)u_i  -q_{ij}u_ix_j-q_{ij}\gamma(x_i)u_j+\theta_{ij}x_ix_j+Kx_i +K'x_j+\hat{\theta}_{ij}p_{ij}(x_1,\cdots,x_n)-\bar{p}_{ij}(u_1,\cdots,u_n)+a_{ij}^{(0)}\theta=0
\end{equation}
where $\hat{\theta}_{ij},\theta_{ij},K,K' \in \A$, $K=0$ if $\delta_j= 0$, and $K'=0$ if $\delta_i=0$.

\begin{proof}
Note that $d(1)=d(1\cdot1)=d(1)+d(1)+\sum_s\theta_s=2d(1)+\theta$, so we have that 
\begin{eqnarray*}
d(x_jx_i) &=& d(q_{ij}x_ix_j+p_{ij}(x_1,\cdots,x_n)) \\[0.15cm]
&=& q_{ij}d(x_i)x_j+q_{ij}\gamma(x_i)d(x_j)+q_{ij}\sum_{s}\theta_s a_{s,i}x_ib_{s,j}x_j + d(a_{ij}^{(0)}) + d(a_{ij}^{(1)}x_1) + \cdots +d(a_{ij}^{(n)}x_n) \\[0.15cm]
&=& q_{ij}d(x_i)x_j+q_{ij}\gamma(x_i)d(x_j) + K_1x_ix_j+ \sum_{s}\theta_sa_{s,i}\delta_i(b_{s,i})x_j+\bar{p}_{ij}(u_1,\cdots,u_n)+a_{ij}^{(0)}d(1) \\[0.15cm]
&=&q_{ij}u_ix_j+q_{ij}\gamma(x_i)u_j + K_1x_ix_j+ K_2x_j+ \bar{p}_{ij}(u_1,\cdots,u_n)-a_{ij}^{(0)}\theta 
\end{eqnarray*}
In the other hand 
\begin{eqnarray*}
d(x_jx_i)&=& u_jx_i+\gamma(x_j)u_i + \sum_{s}\theta_sa_{s,j}x_jb_{s,i}x_i \\[0.15cm]  
&=& u_jx_i+\gamma(x_j)u_i + \sum_{s}\theta_sa_{s,j}b_{s,i}x_jx_i + \sum_{s}\theta_sa_{s,j}\delta_j( b_{s,i})x_i  \\[0.15cm]
&=& u_jx_i+\gamma(x_j)u_i + K_3x_i + K_4 (q_{ij}x_ix_j+p_{ij}(x_1,\cdots,x_n)) 
\end{eqnarray*}
so that 
 \begin{equation*}
 u_jx_i+\gamma(x_j)u_i  -q_{ij}u_ix_j-q_{ij}\gamma(x_i)u_j+\theta_{ij}x_ix_j+Kx_i +K'x_j+\hat{\theta}_{ij}p_{ij}(x_1,\cdots,x_n)-\bar{p}_{ij}(u_1,\cdots,u_n)+a_{ij}^{(0)}\theta=0
 \end{equation*}

\end{proof}

\lemma Let $z\in \O$ then $u_i-[ad_{\gamma}z]x_i$ are solutions of the equation (1).
\begin{proof}
\begin{eqnarray*}
&\p{u_j-[ad_{\gamma}z]x_j}x_i+\gamma(x_j)\p{u_i-[ad_{\gamma}z]x_i}  -q_{ij}\p{u_i-[ad_{\gamma}z]x_i}x_j-q_{ij}\gamma(x_i)\p{u_j-[ad_{\gamma}z]x_j}+& \\[0.15cm]
&\theta_{ij}x_ix_j+Kx_i +K'x_j+\hat{\theta}_{ij}p_{ij}(x_1,\cdots,x_n)-\bar{p}_{ij}(\p{u_1-[ad_{\gamma}z]x_1},\cdots,\p{u_n-[ad_{\gamma}z]x_n})+a_{ij}^{(0)}\theta=& \\[0.15cm]
&=\p{ u_jx_i+\gamma(x_j)u_i  -q_{ij}u_ix_j-q_{ij}\gamma(x_i)u_j+\theta_{ij}x_ix_j+Kx_i +K'x_j+\hat{\theta}_{ij}p_{ij}(x_1,\cdots,x_n)+a_{ij}^{(0)}\theta} + & \\[0.15cm]
& + \p{-[ad_{\gamma}z]x_j}x_i- \gamma(x_j)[ad_{\gamma}z]x_i+q_{ij}\p{[ad_{\gamma}z]x_i}x_j+q_{ij}\gamma(x_i)[ad_{\gamma}z]x_j - \sum_{t=1}^{n}a_{ij}^{(t)}\p{u_t-[ad_{\gamma}z]x_t} =& \\[0.15cm]
&=( u_jx_i+\gamma(x_j)u_i  -q_{ij}u_ix_j-q_{ij}\gamma(x_i)u_j+\theta_{ij}x_ix_j+Kx_i & \\[0.15cm]
&+K'x_j+\hat{\theta}_{ij}p_{ij}(x_1,\cdots,x_n)-\bar{p}_{ij}(u_1,\cdots,u_n)+a_{ij}^{(0)}\theta) \, + & \\[0.15cm]
& \p{-[ad_{\gamma}z]x_j}x_i- \gamma(x_j)[ad_{\gamma}z]x_i+q_{ij}\p{[ad_{\gamma}z]x_i}x_j+q_{ij}\gamma(x_i)[ad_{\gamma}z]x_j + \sum_{t=1}^{n}a_{ij}^{(t)}\p{[ad_{\gamma}z]x_t} =& \\[0.15cm]
&= -zx_jx_i+\gamma(x_j)zx_i-\gamma(x_j)zx_i+\gamma(x_j)\gamma(x_i)z+q_{ij}zx_ix_j-q_{ij}\gamma(x_i)zx_j+q_{ij}\gamma(x_i)zx_j-q_{ij}\gamma(x_i)\gamma(x_j)z+& \\[0.15cm]
&+\sum_{t=1}^na_{ij}^{(t)}\p{zx_t-\gamma(x_t)z} =&\\[0.15cm]
&= [-zx_jx_i+q_{ij}zx_ix_j+\sum_{t=1}^na_{ij}^{(t)}zx_t+a_{ij}^{(0)}z  ]+[\gamma(x_j)zx_i-\gamma(x_j)zx_i]+ & \\[0.15cm]
&+ [-q_{ij}\gamma(x_i)zx_j+q_{ij}\gamma(x_i)zx_j]+[\gamma(x_j)\gamma(x_i)z-\gamma(x_j)\gamma(x_i)z ] +& \\[0.15cm]
&+[ \gamma(x_j)\gamma(x_i)z-q_{ij}\gamma(x_ix_j)z-\sum_{t=1}^na_{ij}^{(t)}\gamma(x_t)z-a_{ij}^{(0)}z]=&\\[0.15cm]
&=[z\p{-x_jx_i+q_{ij}x_ix_j+\sum_{t=1}^na_{ij}^{(t)}x_t+a_{ij}^{(0)}}]+ [ \gamma\p{x_jx_i-q_{ij}x_ix_j-\sum_{t=1}^na_{ij}^{(t)}x_i-a_{ij}^{(0)}}z]=0&
\end{eqnarray*}
\end{proof}

\proposition Let $v=ax_1^{m_1}\cdots x_n^{m_n}$ with $a\in A$ and $m=(m_1,\cdots,m_n)\in Z$ such that either $\gamma_1\neq 1$ and $\p{1-\gamma_1\prod_{j=2}^nq_{j1}^{m_j}}\in A^{\ast}$ or $\gamma_1=1$ and exists $m_j\neq 0$ for $j\neq 1$. Then there exists $w\in A^{\ast}$ with
\begin{equation*}
[ad_{\gamma}wvx_1^{-1}]x_1=v
\end{equation*}

\begin{proof}
Take $w=\p{1-\gamma_1\prod_{j=2}^nq_{j1}^{m_j}}^{-1}$ that exists by hypothesis. Note that $\gamma_1w=w\gamma_1$ so 
\begin{eqnarray*}
[ad_{\gamma}wvx_1^{-1}]x_1 &=& wvx_1^{-1}x_1-\gamma_1x_1wvx_1^{-1} \\[0.15cm]
&=& w \p{ v-\gamma_1ax_1x_1^{m_1}\cdots x_n^{m_n}x_1^{-1}} \\[0.15cm]
&=& w\p{v-\gamma_1ax_1^{m_1}x_1x_2^{m_2}\cdots x_n^{m_n}x_1^{-1}}=w\p{v-\gamma_1ax_1^{m_1}q_{21}^{m_2}x_2^{m_2}x_1\cdots x_n^{m_n}x_1^{-1}} \\[0.15cm]
&=& w\p{v-\gamma_1ax_1^{m_1}q_{21}^{m_2}x_2^{m_2}q_{31}^{m_3}x_3\cdots q_{n1}^{m_n}x_n^{m_n}x_1x_1^{-1} }\\[0.15cm]
&=& w\p{ v - \gamma_1\prod_{j=2}^nq_{j1}^{m_j} a x_1^{m_1}\cdots x_n^{m_n}} = w\p{ 1 - \gamma_1\prod_{j=2}^nq_{j1}^{m_j}}v=v
\end{eqnarray*}
\end{proof} 

\corollary Let $v=\sum_t \rho_tX_t\in\O$  such that $\rho_tX_t$ satisfies the hypothesis of the above proposition for all $t$, where $X_t\in Mon(\O)$. Then there exists $w\in \O$ such that
\[
[ad_{\gamma}w]x_1=v
\]
\begin{proof}
 We have that $[ad_{\gamma}w+w']z=(w+w')z-\gamma(z)(w+w')=\p{wz-\gamma(z)w}+\p{w'z-\gamma(z)w }=[ad_{\gamma}w]z+[ad_{\gamma}w']z$ for all $w,w',z\in \O$ so for each $\rho_tX_t$ there exists $w_t\in A$ such that 
\[ 
[ad_{\gamma}w_t\rho_tX_tx_{1}^{-1}]x_1=\rho_tX_t 
\]
 then
 \[
v=\sum_t\rho_t X_t=\sum_t [ad_{\gamma}w_t\rho_tX_tx_{1}^{-1}]x_1 = [ad_{\gamma}\sum_t w_t\rho_tX_tx_{1}^{-1}]x_1=[ad_{\gamma}w]x_1
 \]
\end{proof}

\theorem Let $u_i=d(x_i) \in \O$  for $i=1,\cdots,n$. 
\begin{enumerate}
\item If $\gamma_1\neq 1$ and $\p{\prod_{j=2}^n q_{j1}^{m_j}-\gamma_1q_{1t}}\in A^{\ast}$ for all $(0,m_2,\cdots,m_n)\in Z $. Then there exists $w\in\O$ and $\rho_j\in A$ such that 
\[
u_1=[ad_{\gamma}w]x_1
\]

\[
u_j=\rho_j x_j+[ad_{\gamma}w]x_j
\]
for all $j\neq 1$ where $\p{q_{1t}-q_{1t}\gamma_1}\rho_t=-\theta_{1t}$.
\item If $\gamma_j=1$ for all $i=1,\cdots,n$. Then there exists $w\in\O$  such that
\[
u_j=\lambda_jx_j+[ad w]x_j
\]
for all $j=1,\dots,n$ and some $\lambda_j=\lambda(x_j)\in A$.
\end{enumerate}
\begin{proof}
\begin{enumerate}
\item By the above corollary  there exists $w\in\O$ such that $u_1=[ad_{\gamma}w]x_1$ so we define $\u_t=u_t-[ad_{\gamma}w]x_t$ where they hold the equation (1). In particular for $i=1$ and $t\neq 1$, $p_{1t}(x_1,\cdots,x_n)=0$ and $\u_1=0$ so 
\begin{equation*}
0=\u_tx_1-q_{1t}\gamma_1x_1\u_t+\theta_{1t}x_1x_t+K_tx_1
\end{equation*}
If $\u_t=\sum_{m\in Z}\eta_{m}x_1^{m_1}\cdots x_n^{m_n}\in\O$ we have 

\begin{eqnarray*}
0 &=& \p{\sum_{m\in Z}\eta_m x_1^{m_1}\cdots x_n^{m_n}x_1-\gamma_1q_{1t} \sum_{m\in Z}\eta_m x_1  x_1^{m_1}\cdots x_n^{m_n}} \\[0.15cm]
&+& \theta_{1t}x_1x_t+Kx_t \\[0.15cm]
&=& \sum_{m\in Z}\eta_m\prod_{j=2}^n q_{1j}^{m_j} x_1x_1^{m_1}\cdots x_n^{m_n}-\gamma_1q_{1t} \sum_{m\in Z}\eta_m    x_1x_1^{m_1}\cdots x_n^{m_n} \\[0.15cm]
&+& \theta_{1t}x_1x_t+Kx_t  \\[0.15cm]
&=& \p{\sum_{m\in Z}\eta_m \p{\prod_{j=2}^n q_{1j}^{m_j}-\gamma_1q_{1t}} x_1^{m_1+1}\cdots x_n^{m_n} }\\[0.15cm]
&+& \theta_{1t}x_1x_t+Kx_t
\end{eqnarray*}
As $Mon(\O)$ is a basis of $\O$ and $\prod_{j=2}^n q_{j1}^{m_j}-\gamma_1q_{1t}\in A^{\ast}$ we have that $\eta_m=0$ if $m\neq (0,\cdots,1,\cdots,0)$ where $1$ is in the $t^{th}
$ position, so $\u_t=\eta_t x_t$ and $K=0$. In the other hand  

\begin{equation*}
0=\eta_t x_tx_1-q_{1t}\gamma_1x_1\eta_tx_t+\theta_{1t}x_1x_t=\p{ \p{q_{1t}-q_{1t}\gamma_1}\eta_t+\theta_{1t}}x_1x_t
\end{equation*}
then $ \p{q_{1t}-q_{1t}\gamma_1}\eta_t+\theta_{1t}=0$

\item Put $u_1=u'_1+v_1$ with $v_1\in\O\setminus A((x_1))$. By the corollary (2.5) there exists $w\in\O$ with $v_1=[ad\, w]x_1$, we will denote $\u_i=u_i-[ad\, w]x_i$ so we have $\u_1\in A((x_1))$. Let $t\neq 1$ 

\begin{equation*}
\u_tx_1-q_{1t}x_1\u_t=-x_t\u_1+q_{1t}\u_1x_t-\theta_{1t}x_1x_t-Kx_t\in A((x_1))x_t
\end{equation*}
because if $\u_1=\sum_{t}\rho_tx_1^t$ then $x_t\u_1=\sum_t\rho_tq_{1t}^tx_1^tx_t$.
If $\u_t=\sum_{m\in Z}\eta_mx_1^{m_1}\cdots x_n^{m_n}$ then 
\begin{eqnarray*}
\u_tx_1-x_1q_{1t}\u_t&=&\sum_{m\in Z}\eta_mx_1^{m_1}\cdots x_n^{m_n}x_1-x_1q_{1t}\sum_{m\in Z}\eta_mx_1^{m_1}\cdots x_n^{m_n} \\[0.15cm]
&=& \sum_{m\in Z} \eta_m \p{\prod_{j=2}^nq_{1j}^{m_j}-q_{1t}}x_1 x_1^{m_1}\cdots x_n^{m_n} \\[0.15cm]
&=& x_1\sum_{m\in Z} \eta_m \p{\prod_{j=2}^nq_{1j}^{m_j}-q_{1t}} x_1^{m_1}\cdots x_n^{m_n}
\end{eqnarray*}
multiplying by $x_1^{-1}$ we have 

\[
\sum_{m\in Z} \eta_m \p{\prod_{j=2}^nq_{1j}^{m_j}-q_{1t}} x_1^{m_1}\cdots x_n^{m_n}\in A((x_1))x_t
\]
Since $\prod_{j=2}^nq_{1j}^{m_j}-q_{1t} \in A^{\ast}$ we have $\eta_m=0$ if $m_i\neq 0$ for $i\neq 1,t$ or if $m_t\neq 1$. Then $\u_t=f_t(x_1)x_t$ where $f_t(x_1)\in A((x_1))$.
Take $t=2$, if $f_2(x_1)=\lambda_2+\sum_{t\neq 0}\rho_tx_1^t$ then for all $t$ put $w_t=\p{1-q_{12}^t}^{-1}$, note that $w_tx_i=x_iw_t $ for all $i=1,\cdots,n$ so  

\begin{eqnarray*}
[ad\, w_t\rho_tx_1^t]x_2 &=& w_t\rho_tx_1^tx_2-x_2w_t\rho_tx_1^t \\[0.15cm]
&=& w_t\p{\rho_tx_1^tx_2-x_2\rho_tx_1^t}=w_t\p{\rho_tx_1^tx_2-\rho_tx_2x_1^t-\delta_2(\rho_t)x_1^t} \\[0.15cm]
&=& w_t\p{1-q_{12}^t}\rho_tx_1^tx_2-w_t\delta_2(\rho_t)x_1^t=\rho_tx_1^tx_2-w_t\delta_2(\rho_t)x_1^t
\end{eqnarray*}

Then there exists $w'\in A((x_1))\setminus A$ or $w'=0$ with $\u_2-[ad\, w']x_2=\lambda_2x_2+f'(x_1)$ and $f'(x_1)\in A((x_1))\setminus A$ or $f'(x_1)=0$.

Take $\U_l=\u_l-[ad\, w']x_l=g_l(x_1)x_l+g'_l(x_1)$,  ($g'_l(x_1)\in A((x_1))\setminus A$ or $g'_l(x_1)=0$ ), because if $w'=\sum_t a_tx_1^t$ then 
\begin{equation*}
[ad\, w']x_l = \sum_ta_tx_1^tx_l-\sum_tx_la_tx_1^t= \sum_ta_tx_1^tx_l-\sum_{t}a_tq_{1t}^tx_1^tx_l-\sum_t\delta_l(a_t)x_1^t  
\end{equation*}
and $\U_2=\lambda_2x_2+g'_2(x_1)$.

By the equation (1) for ($t=1$)  if $\U_1=\sum_t a_tx_1^t$ then 
\begin{eqnarray*}
0 &=& \lambda_2x_2x_1+g'_2(x_1)x_1 + x_2\U_1- q_{12}\U_1x_2-q_{12}x_1\lambda_2x_2-q_{12}x_1g'_2(x_1)+\theta_{12}x_1x_2 + Kx_1 \\[0.15cm]
&=& [\lambda_2q_{12}x_1x_2 + \sum_t x_2a_tx_1^t-\sum_tq_{12}a_tx_1^tx_2-q_{12}\lambda_2x_1x_2 +\theta_{12}x_1x_2] \\[0.15cm]
&+& [g'_2(x_1)x_1+q_{12}x_1g'_2(x_1) +Kx_1] \\[0.15cm]
&=& \p{\lambda_2q_{12}x_1x_2 +\sum_t a_tx_2x_1^t-\sum_tq_{12}a_tx_1^tx_2-q_{12}\lambda_2x_1x_2 +\theta_{12}x_1x_2 } \\[0.15cm]
&+& \p{g'_2(x_1)x_1-q_{12}x_1g'_2(x_1) +Kx_1 + \sum_t \delta_2(a_t)x_1^t } \\[0.15cm]
&=& \p{\sum_t a_t\p{q_{12}^t-q_{12}} x_1^tx_2 + \theta_{12}x_1x_2  } \\[0.15cm]
&+& \p{g'_2(x_1)x_1-q_{12}x_1g'_2(x_1) +Kx_1 + \sum_t \delta_2(a_t)x_1^t }
\end{eqnarray*}
Then $t=1$, $\theta_{12}=0$, and $\U_1=ax_1$. For this equations if $g'_2(x_1)\neq  0 $ put $it(g'_2(x_1))=b_v x_1^v$ with $v\neq 0$ then
\begin{equation*}
0=b_vx_1^{v+1}-q_{12}b_vx_1^{v+1} = b_v(1-q_{12})x_1^{v+1}
\end{equation*}
so $b_v=0$ but this is not possible. Thus $g'_2(x_1)=0$.

Take $3\leq t$  and $\U_1=a_1x_1$.

\begin{eqnarray*}
0 &=& \p{g_t(x_1)x_t+g'_t(x_1)}x_1+x_ta_1x_1  -q_{1t}a_1x_1x_t-q_{1t}x_1\p{g_t(x_1)x_t+g'_t(x_1)}+\theta_{1t}x_1x_t+Kx_1 \\[0.15cm]
&=& \p{g_t(x_1)x_tx_1+a_1x_tx_1-q_{1t}a_1x_1x_t -q_{1t}x_1g_t(x_1)x_t +\theta_{1t}x_1x_t} + \\[0.15cm]
&+& \p{g'_t(x_1)x_1+ \delta_t(a_1)x_1- q_{1t}x_1g'_t(x_1)+Kx_1} \\[0.15cm]
&=& \p{g_t(x_1)x_tx_1-q_{1t}x_1g_t(x_1)x_t +\theta_{1t}x_1x_t} + \p{g'_t(x_1)x_1+ \delta_t(a_1)x_1- q_{1t}x_1g'_t(x_1)+Kx_1} \\[0.15cm]
&=& \theta_{1t}x_1x_t + \p{g'_t(x_1)x_1+ \delta_t(a_1)x_1- q_{1t}x_1g'_t(x_1)+Kx_1}
\end{eqnarray*}
Then $ \theta_{1t}=0$ and if  $g'_t(x_1)\neq 0$ we can take $it(g'_t(x_1))=b_vx_1^v$ with $b_v\neq 0 $. Of this equation we can claim  that $b_v(1-q_{1t})=0$ then $b_v=0$ but this is not possible so $g'_t(x_1)=0$.

In the other hand, if $g_t(x_1)=\sum_t c_tx_1^t$ 

\begin{eqnarray*}
0 &=& g_t(x_1)x_t x_2+x_t\lambda_2x_2 -q_{2t}\lambda_2x_2x_t - q_{2t}x_2g_t(x_1)x_t +\theta_{2t}x_2x_t+ \\[0.15cm]
&+& Kx_2 +K'x_t+\hat{\theta}_{2t}p_{2t}(x_1,\cdots,x_n)-\bar{p}_{2t}(\U_1,\cdots,\U_n)+a_{2t}^{(0)}\theta \\[0.15cm]
&=& g_t(x_1)x_t x_2+\p{\lambda_2x_tx_2 -\lambda_2q_{2t}x_2x_t} - q_{2t}x_2g_t(x_1)x_t +\theta_{2t}x_2x_t+ \\[0.15cm]
&+& Kx_2 +K'x_t+\hat{\theta}_{2t}p_{2t}(x_1,\cdots,x_n)-\bar{p}_{2t}(\U_1,\cdots,\U_n) + \delta_t(\lambda_2)x_2+a_{2t}^{(0)}\theta  \\[0.15cm]
&=& g_t(x_1)\p{q_{2t}x_2x_t+p_{2t}(x_1,\cdots,x_n)} - q_{2t}x_2\sum_t c_tx_1^tx_t +\theta_{2t}x_2x_t +\lambda_2p_{2t}(x_1,\cdots,x_n)+ \\[0.15cm]
&+& Kx_2 +K'x_t+\hat{\theta}_{2t}p_{2t}(x_1,\cdots,x_n)-\bar{p}_{2t}(\U_1,\cdots,\U_n) + \delta_t(\lambda_2)x_2+a_{2t}^{(0)}\theta  \\[0.15cm]
&=& q_{2t}\p{\sum_t c_tx_1^tx_2x_t - \sum_t  c_tx_2x_1^tx_t} +\theta_{2t}x_2x_t + g_t(x_1)p_{2t}(x_1,\cdots,x_n) +\sum_t\delta_2(c_t)x_1^tx_t+ \\[0.15cm]
&+& Kx_2 +K'x_t+\p{\hat{\theta}_{2t}+\lambda_2}p_{2t}(x_1,\cdots,x_n)-\bar{p}_{2t}(\U_1,\cdots,\U_n) + \delta_t(\lambda_2)x_2+a_{2t}^{(0)}\theta   \\[0.15cm]
&=& q_{2t}\p{\sum_t c_t(1-q_{12}^t) x_1^tx_2x_t} +\theta_{2t}x_2x_t + g_t(x_1)p_{2t}(x_1,\cdots,x_n) +\sum_t\delta_2(c_t)x_1^tx_t+ \\[0.15cm]
&+& Kx_2 +K'x_t+\p{\hat{\theta}_{2t}+\lambda_2}p_{2t}(x_1,\cdots,x_n)-\bar{p}_{2t}(\U_1,\cdots,\U_n) + \delta_t(\lambda_2)x_2+a_{2t}^{(0)}\theta   \\[0.15cm]
\end{eqnarray*}
Of this last equation we have that $ c_t(1-q_{12}^t)=0$ for $t\neq 0$ so $c_t=0$ since $(1-q_{12}^t)\in A^{\ast}$,  thus $g_t(x_1)=\lambda_t$ of this way $u_t-[ad\, w+w']x_t=\U_t=\lambda_tx_t $.
\end{enumerate}
\end{proof}

\section{Poisson brackets on $\O$}
\definition A Poisson bracket $\br{\;,\;}$ is a $\Z$-bilinear function on $\O$ such that 
\begin{enumerate}
\item $\br{\;,\;}$  is a Lie bracket.
\item $\br{ab,c}=\br{a,c}b+a\br{b,c}$ for all $a,b,c\in\O$.
\item $0=\br{a,b} +\br{b,a}$ for all $a,b\in\O$.
\end{enumerate}

\remark If $\partial $ is a derivation on $\O$ then for all $i=1,\cdots,r $ we have $0=\partial(1)=\partial(x_ix_i^{-1})=\partial(x_i)x_i^{-1}+x_i\partial(x_i^{-1}) $ then $\partial(x_i^{-1})=-x_i^{-1}\partial(x_i)x_i^{-1} $.

\proposition Let be given a Poisson bracket $\br{a,b}$ on an extension  $\O$ of A. Then there exists $\xi,\xi_j\in A$ such that $\br{x_i,x_j}=\xi[x_i,x_j]+\delta_i(\xi)x_i-\delta_i(\xi)x_j-\delta_i(\xi_{j})$ with $\delta_j(\xi_j)=0 $ for all $ j,i=1,\cdots,n$, where $[a,b]:=ab-ba $.

\begin{proof}
By theorem (2.6) for all $a\in\O$ there exists  $\lambda_a(x_i)\in A$ and $w(a)\in \O$ with 
\[
\br{x_i,a}=\lambda_a(x_i)x_i+[ad\, w(a)]x_i
\]

So we have 
\begin{equation*}
0=\br{x_i,x_i}=\lambda_i(x_i)x_i+[ad\, w(x_i)]x_i
\end{equation*}

Put $w(x_i)=\sum_{m=0}^k\eta_{v_m} x_1^{v_{m_1}}\cdots x_n^{v_{m_n}} $, we will see that $w(x_i)\in A((x_i))$ for $i=1,\cdots,r $ or $w(x_i)\in A[x_i]$ for $i=r+1,\cdots,n$. If $w(x_i)=0 $ we have the claimed, in the other hand, 

\begin{eqnarray*}
0 &=& \lambda_i(x_i)x_i-\sum_{m=0}^kx_i \eta_{v_m} x_1^{v_{m_1}}\cdots x_n^{v_{m_n}} + \sum_{m=0}^k \eta_{v_m} x_1^{v_{m_1}}\cdots x_n^{v_{m_n}} x_i \\[0.15cm]
&=& \lambda_i(x_i)x_i- \sum_{m=0}^k \delta_i(\eta_{v_m}) x_1^{v_{m_1}}\cdots x_n^{v_{m_n}}-\sum_{m=0}^k\eta_{v_m} x_i x_1^{v_{m_1}}\cdots x_n^{v_{m_n}} + \sum_{m=0}^k \eta_{v_m} x_1^{v_{m_1}}\cdots x_n^{v_{m_n}}x_i \\[0.15cm]
&=& \lambda_i(x_i)x_i - \sum_{m=0}^k\delta_i(\eta_{v_m}) x_1^{v_{m_1}}\cdots x_n^{v_{m_n}}-\\[0.15cm]
&-&\sum_{m=0}^k\p{\eta_{v_m} \prod_{s<i}q_{si}^{v_{m_s}} x_1^{v_{m_1}}\cdots x_i^{v_{m_i}+1}\cdots x_n^{v_{m_n}}+ p_{v_m,i}(x_1,\cdots,x_n)}+\\[0.15cm]
&+& \sum_{m=0}^k\p{\eta_{v_m} \prod_{s>i}q_{is}^{v_{m_s}} x_1^{v_{m_1}}\cdots x_i^{v_{m_i}+1}\cdots x_n^{v_{m_n}}+ p_{i,v_m}(x_1,\cdots,x_n)} \\[0.15cm]
&=& \sum_{m=0}^k\eta_{v_m} \p{ \prod_{s>i}q_{is}^{v_{m_s}}-\prod_{s<i}q_{si}^{v_{m_s}}} x_1^{v_{m_1}}\cdots x_i^{v_{m_i}+1}\cdots x_n^{v_{m_n}} \\[0.15cm]
&+& \lambda_i(x_i)x_i- \sum_{m=0}^k \delta_i(\eta_{v_m}) x_1^{v_{m_1}}\cdots x_n^{v_{m_n}} +\sum_{m=0}^k  p_{i,v_m}(x_1,\cdots,x_n)-p_{v_m,i}(x_1,\cdots,x_n)\end{eqnarray*}
where $deg\p{p_{v_m,i}},deg(p_{i,v_m})<v_{m_1}+\cdots+ v_{m_n}+1<\mid v_k\mid+1$ for all $m=0,\cdots,k$. By the last equation we have $\p{\prod_{s<i}q_{si}^{v_{k_s}}- \prod_{s>i}q_{is}^{v_{k_s}}}\eta_{v_k}=0 $ and since $\eta_{v_k}\neq 0$ then $v_{k_l}=0$ for $l\neq i$,  because if there exists $v_{k_l}\neq 0$ for some $l\neq i$ then $\prod_{s<i}q_{si}^{v_{m_s}}- \prod_{s>i}q_{is}^{v_{m_s}}\in A^{\ast}$ by definition of $\O$, so this implies that $\eta_{v_k}=0$, but this is not possible. Later $ p_{i,v_k}=p_{v_k,i}=0$ and 

\begin{eqnarray*}
0 &=& \sum_{m=0}^{k-1}\eta_{v_m} \p{ \prod_{s>i}q_{is}^{v_{m_s}}-\prod_{s<i}q_{si}^{v_{m_s}}} x_1^{v_{m_1}}\cdots x_i^{v_{m_i}+1}\cdots x_n^{v_{m_n}} \\[0.15cm]
&+& \lambda_i(x_i)x_i- \sum_{m=0}^{k-1} \delta_i(\eta_{v_m}) x_1^{v_{m_1}}\cdots x_n^{v_{m_n}} +\sum_{m=0}^{k-1}  p_{i,v_m}(x_1,\cdots,x_n)-p_{v_m,i}(x_1,\cdots,x_n)
\end{eqnarray*}
So $\eta_{v_{k-1}}=0 $, if $\eta_{v_{k-1}}\neq 0 $ then $v_{{k-1}_l}=0 $ for $l\neq i$ because $ \p{\prod_{s<i}q_{si}^{{v_{k-1}}_s}- \prod_{s>i}q_{is}^{{v_{k-1}}_s} }\in A^{\ast}$ if there exist $v_{{k-1}_l}\neq 0 $ for some $l\neq i$. By recurrently way we can claim that $\eta_{v_{m}}=0$ if there exists $v_{m_l}\neq 0 $ for some $l\neq i$, then we have the claimed. Since $w(x_i)=\sum_{t\in\Z}\xi_{it}x_i^t$, where $\xi_{it}=0$ if $t<0$ and $i>r$,  then

\begin{equation*}
0 = \lambda_i(x_i)x_i - x_i\sum_{t\in Z}\xi_{it}x_i^t+\sum_{t\in\Z}\xi_{it}x_i^{t+1}=\lambda_i(x_i)x_i-\sum_{t\in\Z}\delta_{i}(\xi_{it})x_i^t
\end{equation*}
then $\delta_i(\xi_{it})=0$ for $t\neq 1$ and $\lambda_i(x_i)=\delta_i(\xi_{i1}) $. Put  $i< j$, $w(x_i)=\sum_{l=0}^k \xi_{iv_{i_{l}}}x_i^{v_{i_{l}}} $, and $w(x_j)=\sum_{l=0}^r\xi_{jv_{j_{l}}}x_j^{v_{j_{l}}}$, so 

\begin{eqnarray*}
0 &=& \br{x_j,x_i}+\br{x_i,x_j}= \lambda_{i}(x_j)x_j + [ad \, w(x_i)]x_j + \lambda_j(x_i)x_i+[ad\, w(x_j)]x_i \\[0.15cm]
&=& \lambda_{i}(x_j)x_j+ \lambda_j(x_i)x_i + \sum_{l=0}^k \xi_{iv_{i_{l}}}x_i^{v_{i_{l}}}x_j-x_j\sum_{l=0}^k \xi_{iv_{i_{l}}}x_i^{v_{i_{l}}} +  \sum_{l=0}^r\xi_{jv_{j_{l}}}x_j^{v_{j_{l}}}x_i-x_i\sum_{l=0}^r\xi_{jv_{j_{l}}}x_j^{v_{j_{l}}} \\[0.15cm]
&=& \lambda_{i}(x_j)x_j+ \lambda_j(x_i)x_i +  \xi_{iv_{i_k}}x_i^{v_{i_k}}x_j-\xi_{iv_{i_k}}x_jx_i^{v_{i_k}} +  \xi_{jv_{j_r}}x_j^{v_{j_r}}x_i-\xi_{jv_{j_r}}x_ix_j^{v_{j_r}}  \\[0.15cm]
&+&   \sum_{l=0}^{k-1} \xi_{iv_{i_{l}}}x_i^{v_{i_{l}}}x_j-\sum_{l=0}^{k-1}\xi_{iv_{i_{l}}}x_jx_i^{v_{i_{l}}} +  \sum_{l=0}^{r-1}\xi_{jv_{j_{l}}}x_j^{v_{j_{l}}}x_i-\sum_{l=0}^{r-1}\xi_{jv_{j_{l}}}x_ix_j^{v_{j_{l}}}  \\[0.15cm]
&-& \sum_{l=0}^{k-1}\delta_j(\xi_{iv_{i_{l}}})x_i^{v_{i_{l}}} - \sum_{l=0}^{r-1}\delta_i(\xi_{jv_{j_{l}}})x_j^{v_{j_{l}}} -\delta_j(\xi_{iv_{i_k}})x_i^{v_{i_k}}-\delta_i(\xi_{jv_{j_r}})x_j^{v_{j_r}}\\[0.15cm]
&=& \lambda_{i}(x_j)x_j+ \lambda_j(x_i)x_i +  \xi_{iv_{i_k}}x_i^{v_{i_k}}x_j-\xi_{iv_{i_k}}q_{ij}^{v_{i_k}}x_i^{v_{i_k}}x_j +  \xi_{jv_{j_r}}q_{ij}^{v_{j_r}}x_ix_j^{v_{j_r}}-\xi_{jv_{j_r}}x_ix_j^{v_{j_r}}  \\[0.15cm]
&+&   \sum_{l=0}^{k-1} \xi_{iv_{i_{l}}}x_i^{v_{i_{l}}}x_j-\sum_{l=0}^{k-1}\xi_{iv_{i_{l}}}q_{ij}^{v_{i_{l}}}x_i^{v_{i_{l}}}x_j +  \sum_{l=0}^{r-1}\xi_{jv_{j_{l}}}q_{ij}^{v_{j_{l}}}x_ix_j^{v_{j_{l}}}-\sum_{l=0}^{r-1}\xi_{jv_{j_{l}}}x_ix_j^{v_{j_{l}}}  \\[0.15cm]
&-& \sum_{l=0}^{k}\delta_j(\xi_{iv_{i_{l}}})x_i^{v_{i_{l}}} - \sum_{l=0}^{r}\delta_i(\xi_{jv_{j_{l}}})x_j^{v_{j_{l}}} - \xi_{iv_{i_k}}p_{v_{i_k},j}+\xi_{jv_{j_r}}p_{i,v_{j_r}} - \sum_{l=0}^{k-1}\xi_{iv_{i_{l}}} p_{v_{i_{l}},j}+\sum_{l=0}^{r-1}\xi_{jv_{j_{l}}}p_{i,v_{j_{l}}}  \\[0.15cm]
&=& \lambda_{i}(x_j)x_j+ \lambda_j(x_i)x_i +  \xi_{iv_{i_k}}\p{1-q_{ij}^{v_{i_k}}}x_i^{v_{i_k}}x_j + \xi_{jv_{j_r}}\p{q_{ij}^{v_{j_r}}-1}x_ix_j^{v_{j_r}}  \\[0.15cm]
&+&   \sum_{l=0}^{k-1} \xi_{iv_{i_{l}}}\p{1-q_{ij}^{v_{i_{l}}}}x_i^{v_{i_{l}}}x_j +   \sum_{l=0}^{r-1}\xi_{jv_{j_{l}}}\p{q_{ij}^{v_{j_{l}}}-1}x_ix_j^{v_{j_{l}}} -  \\[0.15cm]
&-& \sum_{l=0}^k\delta_j(\xi_{iv_{i_{l}}})x_i^{v_{i_{l}}} - \sum_{l=0}^r\delta_i(\xi_{jv_{j_{l}}})x_j^{v_{j_{l}}} -\sum_{l=0}^{k}\xi_{iv_{i_{l}}} p_{v_{i_{l}},j}+\sum_{l=0}^{r}\xi_{jv_{j_{l}}}p_{i,v_{j_{l}}}
\end{eqnarray*}
where $deg(p_{v_{i_l},j})<v_{i_l}+1<v_{i_k}+1$ and $deg( p_{i,v_{j_l}})<v_{j_l}+1<v_{j_r}+1 $ for  all $l$. Since $0=\xi_{iv_{i_k}}\p{1-q_{ij}^{v_{i_k}}}x_i^{v_{i_k}}x_j + \xi_{jv_{j_r}}\p{q_{ij}^{v_{j_r}}-1}x_ix_j^{v_{j_r}} $, we have $ v_{j_r}=v_{i_k}$ and $v_{j_r}\in \{0,1\} $, later $p_{v_{i_k},j},\; p_{i,v_{j_r}} \in A+Ax_1+\cdots+Ax_n   $.  
The same way
\begin{eqnarray*}
0 &=&   \lambda_{i}(x_j)x_j+ \lambda_j(x_i)x_i +  \xi_{iv_{i_{k-1}}}\p{1-q_{ij}^{v_{i_{k-1}}}}x_i^{v_{i_{k-1}}}x_j + \xi_{jv_{j_{r-1}}}\p{q_{ij}^{v_{j_{r-1}}}-1}x_ix_j^{v_{j_{r-1}}}  \\[0.15cm]
&+&   \sum_{l=0}^{k-2} \xi_{iv_{i_{l}}}\p{1-q_{ij}^{v_{i_{l}}}}x_i^{v_{i_{l}}}x_j +   \sum_{l=0}^{r-2}\xi_{jv_{j_{l}}}\p{q_{ij}^{v_{j_{l}}}-1}x_ix_j^{v_{j_{l}}} - \\[0.15cm]
&-& \sum_{l=0}^{k}\delta_j(\xi_{iv_{i_{l}}})x_i^{v_{i_{l}}} - \sum_{l=0}^{r}\delta_i(\xi_{jv_{j_{l}}})x_j^{v_{j_{l}}} - \sum_{l=0}^{k}\xi_{iv_{i_{l}}} p_{v_{i_{l}},j}+\sum_{l=0}^{r}\xi_{jv_{j_{l}}}p_{i,v_{j_{l}}}
 \end{eqnarray*} 
Then $0=\xi_{iv_{i_{k-1}}}\p{1-q_{ij}^{v_{i_{k-1}}}}x_i^{v_{i_{k-1}}}x_j + \xi_{jv_{j_{r-1}}}\p{q_{ij}^{v_{j_{r-1}}}-1}x_ix_j^{v_{j_{r-1}}}$, since  $v_{i_{k-1}},\;v_{j_{r-1}}\leq 0$, $\xi_{jv_{j_{k-1}}},\;\xi_{iv_{i_{r-1}}} \neq 0$, and $\p{1-q_{ij}^{v_{i_{k}}}}, \p{q_{ij}^{v_{j_{k}}}-1}\in A^{\ast}$ if $v_{i_{k-1}},\; v_{j_{r-1}}\neq 0$ then we can claim that $v_{i_{k-1}},\; v_{j_{r-1}}= 0$. Of recurrently form we can claim that $w(x_t)=\xi_{t1}x_t+\xi_{t0}$ with $\xi_{t1},\xi_{t0}\in A $ and $\delta_t(\xi_{t0})=0 $ for all $t=1,\cdots,n$. Again put $ i<j$
\begin{eqnarray*}
0 &=& \br{x_j,x_i}+\br{x_i,x_j}= \lambda_{i}(x_j)x_j + [ad \, w(x_i)]x_j + \lambda_j(x_i)x_i+[ad\, w(x_j)]x_i \\[0.15cm]
&=& \lambda_{i}(x_j)x_j+ \lambda_j(x_i)x_i + \xi_{i1}x_ix_j+\xi_{i0}x_j-x_j\xi_{i1}x_i - x_j\xi_{i0} \\[0.15cm] 
&+& \xi_{j1}x_jx_i+\xi_{j0}x_i-x_i\xi_{j1}x_j - x_i\xi_{j0} \\[0.15cm]
&=& \lambda_{i}(x_j)x_j+ \lambda_j(x_i)x_i + \xi_{i1}x_ix_j+\xi_{i0}x_j- \p{\delta_{j}(\xi_{i1})x_i +\xi_{i1}q_{ij}x_ix_j+\xi_{i1}p_{ij}(x_1,\cdots,x_n) }- \\[0.15cm]
&-&\p{\xi_{i0}x_j+\delta_j(\xi_{i0})}+\p{\xi_{j1}q_{ij}x_ix_j+\xi_{j1}p_{ij}(x_1,\cdots,x_n) } + \\[0.15cm]
&+& \xi_{j0}x_i-\p{\delta_i(\xi_{j1})x_j+ \xi_{j1}x_ix_j  }-\p{\delta_i(\xi_{j0})+\xi_{j0}x_i  }  \\[0.15cm]
&=&  \xi_{i1}x_ix_j - \xi_{i1}q_{ij}x_ix_j  + \xi_{j1}q_{ij}x_ix_j -\xi_{j1}x_ix_j \\[0.15cm]
&+& \lambda_{i}(x_j)x_j+ \lambda_j(x_i)x_i  -\delta_j(\xi_{i1})x_i-\delta_{j}(\xi_{i0})-\xi_{i1}p_{ij}(x_1,\cdots,x_n) \\[0.15cm]
&+& \xi_{j1}p_{ij}(x_1,\cdots,x_n)-\delta_i(\xi_{j1})x_j-\delta_i(\xi_{j0})
\end{eqnarray*}
So $(\xi_{i1}-\xi_{j1})(1-q_{ij})=0$ as $1-q_{ij}\in A^{\ast} $ then $\xi_{1j}=\xi_{1i}=\xi $. Later 
\begin{equation*}
0=(\lambda_{i}(x_j)-\delta_j(\xi))x_j+ \p{\lambda_j(x_i)-\delta_i(\xi)}(x_i)x_i -\delta_j(\xi_{i0})-\delta_i(\xi_{j0})
\end{equation*}
And $\lambda_{i}(x_j) = \delta_j(\xi) $, $\lambda_j(x_i)=\delta_i(\xi)$, $  \delta_j(\xi_{i0})=-\delta_i(\xi_{j0})$. Note that 
\begin{eqnarray*}
\br{x_i,x_j}&=&\lambda_{j}(x_i)x_i+[ad\, \xi x_j+\xi_{j0}]x_i=\lambda_{j}(x_i)x_i+\xi x_jx_i+\xi_{j0}x_i-x_i\xi x_j-x_i\xi_{j0} \\[0.15cm]
&=& \delta_i(\xi)x_i+\xi x_jx_i+\xi_{j0}x_i - \xi x_ix_j + \delta_i(\xi)x_j-\xi_{j0}x_i-\delta_i(\xi_{j0}) \\[0.15cm]
&=& \xi[x_i,x_j]+\delta_i(\xi)x_i-\delta_i(\xi)x_j-\delta_i(\xi_{j0})
\end{eqnarray*}
\end{proof}

\lemma Let $\br{\;,\;}$ be a Poisson bracket on $\O$ where $\delta_i=0$ for all $i=1,\cdots,n$.  Then there exists $\xi\in A$ such that $\br{x_i,a}=\xi[x_i,a]$ for all $i=1,\cdots,n$ and $ a=x_{t_1}^{b_1}\cdots x_{t_n}^{b_n}\in \O$ monomial.
\begin{proof}
Let $\xi\in A $ such that $\br{x_i,x_j}=\xi[x_i,x_j]$  for all $i,j=1,\cdots,n$ which is given by proposition (3.3) and $a=x_{t_1}^{b_1}\cdots x_{t_l}^{b_l}$ with $b_r\neq 0$ for all $r=1,\cdots,l$. Fix $i=1,\cdots,n$, we will show, by induction on $l$, that $\br{a,x_i}=\xi[a,x_i]$.

\begin{enumerate}
\item  ($a=x_t^n$) For this case we will do induction on $n$.
\begin{enumerate}
\item ($a=x_t$) we have $\br{x_t,x_i}=\xi[x_t,x_i]$ by the above proposition.
\item ($a=x_t^{n+1}$) 
\begin{eqnarray*}
\br{x_t^{n+1},x_i}&=&\br{x_t^n,x_i}x_t+x_t^n\br{x_t,x_i} \\[0.15cm]
&=& \p{\xi[x_t^n, x_i]}x_t+x_t^n\p{\xi[ x_t,x_i]} \\[0.15cm]
&=& \xi\p{-x_ix_t^{n+1}+x_t^nx_ix_t-x_t^nx_ix_t+x_t^{n+1}x_i  }=\xi[x_t^{n+1},x_i]
\end{eqnarray*}
\item ($a=x_t^{-1}$, if $t\leq r$)
\begin{eqnarray*}
\br{x_t^{-1},x_i}&=& -x_t^{-1}\br{x_t,x_i}x_t^{-1} \\[0.15cm]
&=& -\xi x_t^{-1}\p{[x_t,x_i]} x_t^{-1} =\xi\p{-x_t^{-1}\p{x_tx_i-x_ix_t}x_t^{-1}  } \\[0.15cm]
&=& \xi\p{x_t^{-1}x_i-x_ix_t^{{-1}} }=\xi[x_t^{-1}, x_i]
\end{eqnarray*}
\item ($a=x_t^{-n-1}$, if $t\leq r$) 
\begin{eqnarray*}
\br{x_t^{-1-n},x_i}&=&\br{x_t^{-n},x_i}x_t^{-1}+x_t^{-n}\br{x_t^{-1},x_i} = \xi\p{ \p{[x_t^{-n},x_i]}x_t^{-1} + x_t^{-n}\p{[x_t^{-1},x_i]}} \\[0.15cm]
&=& \xi\p{-x_ix_t^{-n-1}+x_t^{-n}x_ix_t^{-1}-x_t^{-n}x_ix_t^{-1}+x_t^{-1-n}x_i}=\xi[x_t^{-n-1},x_i]
\end{eqnarray*}
\end{enumerate} 

\item ($a=x_{t_1}^{b_1}\cdots x_{t_l}^{b_l}$) 
\begin{eqnarray*}
\br{a,x_i} &=& \br{x_{t_1}^{b_1}\cdots x_{t_l}^{b_l},x_i} = \br{x_{t_1}^{b_1},x_i}x_{t_2}^{b_2}\cdots x_{t_l}^{b_l}+ x_{t_1}^{b_1}\br{x_{t_2}^{b_2}\cdots x_{t_l}^{b_l},x_i} \\[0.15cm]
&=& \xi\p{\p{[x_{t_1}^{b_1}, x_i]}x_{t_2}^{b_2}\cdots x_{t_l}^{b_l}+x_{t_1}^{b_1}\p{[x_{t_2}^{b_2}\cdots x_{t_l}^{b_l}, x_i]}} \\[0.15cm]
&=& \xi\p{-x_ia+x_{t_1}^{b_1}x_ix_{t_2}^{b_2}\cdots x_{t_l}^{b_l} - x_{t_1}^{b_1}x_ix_{t_2}^{b_2}\cdots x_{t_l}^{b_l} + ax_i }=\xi[a, x_i]
\end{eqnarray*}
\end{enumerate}
So we have the claimed. By theorem (2.6) 

\begin{eqnarray*}
0 &=& \br{x_i,a}+\br{a,x_i}=\lambda_a(x_i)x_i+[ad\, w(a)]x_i+\br{a,x_i} \\[0.15cm]
&=& \lambda_a(x_i)x_i + [ad\, w(a)]x_i - \xi[ad\, a]x_i = \lambda_a(x_i)x_i + [ad\, w(a)-\xi a]x_i 
\end{eqnarray*}
but if we take $w(a)-\xi a$ as $ w(x_i)$ in the proof of the above proposition,  we can show that $w(a)-\xi a \in A((x_i)) $ ( if $ i=1,\cdots,r$ ) or  $w(a)-\xi a \in A[x_i] $ ( if $i=r+1,\cdots,n $) for all $i=1,\cdots,n$ then $ w(a)-\xi a =l$ with $ l\in A$. Then $[ad\, w(a)]x_i=[ad\, \xi a+l]x_i=[ad \, \xi a]x_i+[ad\, l]x_i=\p{\xi a x_i-x_i\xi a}+\p{ lx_i-x_il } =\xi[x_i,a]$ and $\lambda_a(x_i)=0$ for all $i=1\cdots,n$.
\end{proof}

\theorem Let $\br{\,,\,}$ be  a Poisson bracket on $\O$. If $A$ is a commutative ring,  $\delta_i=0$ for all $i=1,\cdots,n$,   and $\br{\,,\,}$ is a A-bilinear function, then there exists $\xi\in A$ such that $ \br{a,b}=\xi[a,b ]$ for all $a,b \in\O$.

\begin{proof}
Let $a,b\in\O$ and $\xi\in A$ such that $\br{x_i,X}=\xi[x_i,X]$ for all $i=1,\cdots,n$, where $X\in Mon\{\O\}$, $\xi$ is given by  lemma (3.4). If $b=\sum_{r=0}^{k}\eta_{v_r}x_{1}^{v_{r_1}}\cdots x_n^{v_{r_n}} $  
since $\br{\,,\,}$ is a $A$-bilinear function and $x_ic=cx_i+\delta_i(c)=cx_i$ for all $c\in A$ and $ i=1,\cdots,n$,  it is enough to see this when $a =x_{t_1}^{a_1}\cdots x_{t_l}^{a_l} $ with $a_t\neq 0$. We will show this by induction on $l$. 
\begin{enumerate}
\item ($a=x_t^n$) For this case we will do induction on $n$.
\begin{enumerate}
\item ($a=x_t$) By the above lemma we have that 
\begin{eqnarray*}
\br{x_{t},b}&=&\br{x_t,\sum_{r=0}^{k}\eta_{v_r}x_{1}^{v_{r_1}}\cdots x_n^{v_{r_n}}}=\sum_{r=0}^{k}\eta_{v_r}\br{x_t,x_{1}^{v_{r_1}}\cdots x_n^{v_{r_n}}}\\[0.15cm]
&=&\sum_{r=0}^{k}\eta_{v_r}\xi[x_t,x_{1}^{v_{r_1}}\cdots x_n^{v_{r_n}}]=\xi\sum_{r=0}^{k}\eta_{v_r}[x_t,x_{1}^{v_{r_1}}\cdots x_n^{v_{r_n}}]\\[0.15cm]
&=& \xi[x_t,b]
\end{eqnarray*}
\item ($a=x_t^{n+1}$) 
\begin{eqnarray*}
\br{x_t^{n+1},b}&=&\br{x_t^n,b}x_t+x_t^n\br{x_t,b} = \xi\p{[x_t^n,b]x_t + x_t^n[x_t,b]} \\[0.15cm]
&=& \xi\p{x_t^nbx_t-bx_t^{n+1}+ x_t^{n+1}b- x_t^nbx_t } =\xi[x_t^{n+1},b]
\end{eqnarray*}

\item ($a=x_t^{-1}$) 
\[\br{x_t^{-1},b}=-x_t^{-1}\br{x_t,b}x_t^{-1}=\xi\p{-x_t^{-1}[x_t,b]x_t^{-1}}=\xi[x_t^{-1},b]
\]

\item ($a=x_t^{-1-n}$) 
\begin{eqnarray*}
\br{x_t^{-1-n},b} &=& \br{x_t^{-n},b}x_t^{-1}+x_t^{-n}\br{x_t^{-1},b}=\xi\p{[x_t^{-n},b]x_t^{-1}+ x_t^{-n}[x_t^{-1},b]} \\[0.15cm]
&=& \xi\p{x_t^{-n}bx_t^{-1}-bx_t^{-n-1}+x_t^{-n-1}b-x_t^{-n}bx_t^{-1}}=\xi[x_t^{-n-1},b]
\end{eqnarray*}
\end{enumerate}

\item ($a=x_{t_1}^{a_1}\cdots x_{t_l}^{a_l}$)
\begin{eqnarray*}
\br{a,b}&=&\br{x_{t_1}^{a_1}\cdots x_{t_l}^{a_l},b} = \br{x_{t_1}^{a_1},b}x_{t_2}^{a_2}\cdots x_{t_l}^{a_l}+x_{t_1}^{a_1}\br{x_{t_2}^{a_2}\cdots x_{t_l}^{a_l},b} \\[0.15cm]
&=& \xi\p{[x_{t_1}^{a_1},b]x_{t_2}^{a_2}\cdots x_{t_l}^{a_l}+x_{t_1}^{a_1}[x_{t_2}^{a_2}\cdots x_{t_l}^{a_l},b] } \\[0.15cm]
&=& \xi\p{x_{t_1}^{a_1}bx_{t_2}^{a_2}\cdots x_{t_l}^{a_l}-ba+ab-x_{t_1}^{a_1}bx_{t_2}^{a_2}\cdots x_{t_l}^{a_l}  } =\xi[a,b]
\end{eqnarray*}

\end{enumerate} 
\end{proof}

\lemma Let $B=\sigma(A)\langle x_1,\cdots,x_n\rangle$ be a skew PBW extension 
 of $A$ with 

\begin{enumerate}
\item $\sigma_i $ is the identity of $A$ for all $i=1,\cdots,n$.
\item $q_{ij},a_{ij}^{(t)}\in Z(A)$ for all $1\leq i,j \leq n$ and $t=0,\cdots,n$.
\item $\delta_1= 0$
\item $\delta_t(q_{ij})=\delta_t(a_{ij}^{(m)}) =0$ for all $m = 0, \cdots, n$ and  $i,j,t = 1, \cdots, n$.
\item $p_{1j}=0$ for all $j=1,\cdots ,n$.
\end{enumerate}  then there exists $\o^{1,n}$ with the same properties.
\begin{proof}
We will see  that $BS^{-1}$  exists, showing the Ore conditions on $S$, and it is an extension of $A$ of type $\o^{1,n}$
for some set $S$. Let $S:=\{x_1^{m}\vert m\in\N\}$, S is a multiplicative set of $B$. We will see that $S$ is a right (left) Ore set.
\begin{enumerate}
\item (Right) Take $P\in B$ and $s\in S$ with $sP=0$, if $P=\sum_{t\in Z}b_t x_1^{t_1}\cdots x_n^{t_n}$ and $s=x_1^l $ we have $ 0= \sum_{t\in Z}b_tx_1^lx_1^{t_1}\cdots x_n^{t_n}$ then $b_t=0$ for all $t\in Z$ then $P=0$ and $Ps=0$.
\item (Left) If $Ps=0$, $0=\sum_{t\in Z}b_t x_1^{t_1}\cdots x_n^{t_n} x_1^l=\sum_{t\in Z}b_t \prod_{j=2}^n q_{1j}^{t_j}x_1^{t_1+l}x_2^{t_2}\cdots x_n^{t_n} $ since $q_{ij}\in A^{\ast} $  then $b_t=0 $ for all $t\in Z$ and $sP=0$.

\item (Right) Let $P=\sum_{t\in Z}b_t x_1^{t_1}\cdots x_n^{t_n} \in B$ and $x_1^l\in S$ then 
\[P x_1^l = \sum_{t\in Z}b_t \prod_{j=2}^n q_{1j}^{t_j}x_1^{t_1+l}x_2^{t_2}\cdots x_n^{t_n} =x_1^l\p{\sum_{t\in Z}b_t \prod_{j=2}^n q_{1j}^{t_j}x_1^{t_1}x_2^{t_2}\cdots x_n^{t_n}}
\]

\item (Left) Let $P=\sum_{t\in Z}b_t x_1^{t_1}\cdots x_n^{t_n} \in B$ and $x_1^l\in S$ then 

\[
x_1^lP=\sum_{t\in Z}b_t x_1^lx_1^{t_1+l}x_2^{t_2}\cdots x_n^{t_n} = \p{ \sum_{t\in Z}b_t \prod_{j=2}^n q_{j1}^{t_j}x_1^{t_1}x_2^{t_2}\cdots x_n^{t_n}} x_1^l
\]

Then $BS^{-1}$ and $S^{-1}B$ exist, so $ BS^{-1}\cong S^{-1}B$. We will see that $D=BS^{-1}$ is an extension of A of type $\o^{1,n}$ and holds the conditions.    

\end{enumerate}
\begin{enumerate}
\item $A \hookrightarrow B  \hookrightarrow D$
\item $\dfrac{x_i}{1}\dfrac{a}{1}=\dfrac{x_ia}{1}=\dfrac{ax_i+\delta_i(a)}{1}=\dfrac{ax_i}{1}+\dfrac{\delta_i(a)}{1} = \dfrac{a}{1}\dfrac{x_i}{1}+\dfrac{\delta_i(a)}{1} $ for $a\in A$.
\item Let $i<j$

\[\dfrac{x_j}{1}\dfrac{x_i}{1}=\dfrac{x_jx_i}{1}=\dfrac{q_{ij}x_ix_j+a_{ij}^{(0)}+a_{ij}^{(1)}x_1+\cdots+a_{ij}^{(n)}x_n}{1}=\dfrac{q_{ij}}{1}\dfrac{x_i}{1}\dfrac{x_j}{1}+\dfrac{a_{ij}^{(0)}}{1}+\dfrac{a_{ij}^{(1)}}{1}\dfrac{x_1}{1}+\cdots +\dfrac{a_{ij}^{(n)}}{1}\dfrac{x_n}{1}
\]
\item Let $\dfrac{\sum_{t\in Z}b_t x_1^{t_1}\cdots x_n^{t_n}}{x_1^l}\in B$, note that $\dfrac{1}{x_1}\dfrac{a}{1}=\dfrac{a}{1}\dfrac{1}{x_1} \;\;$ $\p{\dfrac{1}{x_1} \p{\dfrac{x_1}{1}\dfrac{a}{1}} \dfrac{1}{x_1}= \dfrac{1}{x_1} \p{\dfrac{a}{1}\dfrac{x_1}{1}}  \dfrac{1}{x_1}}$  then 
\begin{equation*}
\dfrac{\sum_{t\in Z}b_t x_1^{t_1}\cdots x_n^{t_n}}{x_1^l} = \p{\dfrac{x_1}{1}}^{-l}\dfrac{\sum_{t\in Z}b_t x_1^{t_1}\cdots x_n^{t_n}}{1}=\sum _{t\in Z} \dfrac{b_t}{1}\p{\dfrac{x_1}{1}}^{t_1-l}\cdots\p{\dfrac{x_n}{1}}^{t_n}
\end{equation*}
\item  If  $P=\sum _{t\in Z} \dfrac{b_t}{1}\p{\dfrac{x_1}{1}}^{t_1-l}\cdots\p{\dfrac{x_n}{1}}^{t_n}=\dfrac{0}{1} $  then  $\dfrac{0}{1}=\p{\dfrac{x_1}{1}}^lP=\dfrac{\sum_{t\in Z}b_t x_1^{t_1}\cdots x_n^{t_n}}{1}$ then there exists $x_1^l\in S $ with 
\[ 0= \sum_{t\in Z}b_t x_1^{t_1}\cdots x_n^{t_n} x_1^l=\sum_{t\in Z}b_t\prod_{j=2}^nq_{1j}^{t_j} x_1^{t_1+l}\cdots x_n^{t_n}
\]
As $ \prod_{j=2}^nq_{1j}^{t_j} \in A^{\ast}$ then $ b_t=0$ for all $t\in Z$ so if we denote $\dfrac{x_i}{1}:=x_i$, then $Mon\{x_1^{\pm 1},x_2,\cdots,x_n\}$ is a $A$-basis of $D$. Note that if we put $\bar{\delta}_i(\dfrac{a}{1})=\dfrac{\delta_i(a)}{1}$ then $\bar{\delta}_i\p{\dfrac{a_{ij}^{(t)}}{1}},\;\bar{\delta}_i\p{\dfrac{q_{ij}}{1}  }=0$ for all $i,j=1,\cdots,n$ and $t=0,\cdots,n$,  and $\bar{\delta}_1=0$.
\end{enumerate}
\end{proof}
\lemma Let $\br{\,,\,}$ be a Poisson bracket over a ring A then $\br{a^l,a^r}=0$ for all $a\in A$ and $l,r\in \N$. 
\begin{proof}
Take $l=1$ and $r=1$ then $\br{a,a}=0$, now $\br{a,a^{l+1}}=\br{a,a}a^l+a\br{a,a^l}=0$ by induction, later $\br{a^r,a^{l+1}}=\br{a^r,a^l}a+a^l\br{a^r,a}=-a^l\br{a,a^r}=0$
\end{proof}

\proposition  Let $B=\sigma(A)\langle x_1,\cdots,x_n\rangle$ be a skew PBW extension of a commutative ring  $A$ such that 
\begin{enumerate}
\item $\sigma_i $ is the identity of $A$ for all $i=1,\cdots$,n.
\item $\delta_i= 0$ for all $i=1,\cdots,n$.
\item For every $i$ fixed, with $i=0,\cdots,n$ and $(m_1,\cdots,m_n)\in Z\setminus \{(0,\cdots,0\}$, $\p{1-\prod_{j=1,j\neq i}^n q_{ij}^{m_j} }\in \A^{\ast}$.
\item $p_{1j}=0$ for all $j=1,\cdots ,n$.
\end{enumerate}
and $\br{\,,\,}$ a Poisson bracket  on $B$ then there exists a Poisson bracket  $\br{\,,\,}_0$ on $\o^{1,n}$  such that $\br{\,,\,}_0\vert_{B}=\br{\,,\,} $. 
\begin{proof}
Let $D=\o^{1,n}$, which is given by lemma (3.6). Put $s\in B$ and define $g_s(x_t)=\br{x_t,s} $ for $t=1,\cdots,n$,   $g_s(x_1^{-1})=-x_1^{-1}\br{x_1,s}x_1^{-1} $, and
\[
g_s(x_{t_1}^{a_1}\cdots x_{t_n}^{a_n})=g_s(x_{t_1}^{a_1})x_{t_2}^{a_2}\cdots x_{t_n}^{a_n}+ x_{t_1}^{a_1}g_s(x_{t_2}^{a_2}\cdots x_{t_n}^{a_n})
\]
of recurrently form.

By the universal property of basis  there exists an A-homomorphism on  $D$ to itself, $g_s$ with $ g_s\vert_{B}=\br{\;,s}$. We will prove that this is a derivation.  
Since a $g_s$ is an A-homomorphism then it is enough to see this to products of monomials. Let $p=x_{t_1}^{a_1}\cdots x_{t_s}^{a_s}$ and $q=x_{l_1}^{b_1}\cdots x_{l_r}^{b_r} $, we see the claimed by induction on $s$.

\begin{enumerate}
\item Let $s=1$, we will do induction on $r$
\begin{enumerate}
\item ($r=1$) We will denote $p=x_i^a$ and $q=x_j^b$.
\begin{enumerate}
\item  ($i<j$) We have 
\[
g_s(x_i^ax_j^b):=g_s(x_i^a)x_j^b+x_i^ag_s(x_j^b)
\]
\item ($x_i^{a},x_j^b\in B $) Then 
\[
g_s(x_i^ax_j^b):=\br{x_i^ax_j^b,s}=\br{x_i^a,s}x_j^b+x_i^a\br{x_j^b,s} =g_s(x_i^a)x_j^b+x_i^ag_s(x_j^b)
\]

\item ($x_1^b\not\in B$, $x_j^a\in B$, $1\leq j$) Note that for all $r\in A$ we have 
\[x_1^{-1}r=x_1^{-1}\p{rx_1 }x_1^{-1}=
 x_1^{-1}\p{x_1r }x_1^{-1}=rx_1^{-1}
\]
and  
\[x_jx_1^{-1} =x_1^{-1}\p{x_1x_j }x_1^{-1}=x_1^{-1}\p{q_{j1}x_jx_1 }x_1^{-1}=q_{j1}x_jx_1^{-1}
 \]
then for all $b\in\Z^-$ and $a\in \N$, $x_j^ax_1^{b}=q_{j1}^{-ab}x_1^{ b}x_j^a=c^{-1}x_1^{b}x_j^a$ where $c\in Z(B)$ and 
\[
c^{-1}g_s(x_j^{a})x_1^{-b}+c^{-1}x_j^{a}g_s(x_1^{-b})=g_s(c^{-1}x_j^ax_1^{-b})=g_s(x_1^{-b}x_j^a)=g_s(x_1^{-b})x_j^a+x_1^{-b}g_s(x_j^a)
\]
\begin{eqnarray*}
\Pa(x_j^ax_1^b)&=&\Pa(c^{-1}x_1^{b}x_j^{a})\\[0.15cm]
&=& c^{-1}\p{\Pa(x_1^b)x_j^a+x_1^b\Pa(x_j^a)}\\[0.15cm]
&=& c^{-1}\p{-x_1^b\pa(x_1^{-b})x_1^bx_j^a+x_1^b\pa(x_j^a)} \\[0.15cm]
&=&-x_1^b\pa(x_1^{-b})c^{-1}x_1^bx_j^a+c^{-1}x_1^b\pa(x_j^a)\\[0.15cm]
&=& -x_1^b\pa(x_1^{-b})x_j^ax_1^b+c^{-1}x_1^b\pa(x_j^a)\\[0.15cm]
&=&x_1^b\p{-\pa(x_1^{-b})x_j^a + c^{-1}\pa(x_j^a)x_1^{-b}}x_1^b \\[0.15cm]
&=&x_1^b\p{x_1^{-b}\pa(x_j^a)-c^{-1}x_j^a\pa(x_1^{-b})}x_1^b \\[0.15cm]
&=&\pa(x_j^a)x_1^b-c^{-1}x_1^bx_j^a\pa(x_1^{-b})x_1^b \\[0.15cm]
&=&\pa(x_j^a)x_1^b-x_j^ax_1^b\pa(x_1^{-b})x_1^b \\[0.15cm]
&=& \Pa(x_j^a)x_1^b+x_j^a\Pa(x_1^{b}) 
\end{eqnarray*}

\end{enumerate} 

\item Let $q= x_1^{l_1}x_2^{l_2}\cdots x_n^{l_n}$ and $p=x_1^b$  then 

\begin{eqnarray*}
\Pa(pq) &=& \pa(x_1^{l_1+b}x_2^{l_2}\cdots x_n^{l_n}) \\[0.15cm]
&:=& \pa(x_1^{l_1+b})x_2^{l_2}\cdots x_n^{l_n}+x_1^{l_1+b}\pa(x_2^{l_2}\cdots x_n^{l_n}) \\[0.15cm]
&=&  \p{\pa(x_1^b)x_1^{l_1}+x_1^b\pa(x_1^{l_1})}x_2^{l_2}\cdots x_n^{l_n}+x_1^{l_1+b}\pa(x_2^{l_2}\cdots x_n^{l_n}) \\[0.15cm]
&=& \pa(p)q+x_1^b\p{\pa(x_1^{l_1})x_2^{l_2}\cdots x_n^{l_n}+x_1^{l_1}\pa(x_2^{l_2}\cdots x_n^{l_n})} \\[0.15cm]
&=& \pa(p)q+p\pa(q)
\end{eqnarray*}  
  
\end{enumerate}
 
\item ($q= x_1^{l_1}x_2^{l_2}\cdots x_n^{l_n}$, $p=x_{t_1}^{a_1}\cdots x_{t_s}^{a_s}$) Let $T\in D$ such that
$x_2^{t_2}\cdots x_n^{t_n}q=T$,   later 
\begin{eqnarray*}
\Pa(pq)&=& \Pa(x_{t_1}^{a_1}x_{t_2}^{a_2}\cdots x_{t_s}^{a_s}q) =\Pa(x_{t_1}^{a_1}T) \\[0.15cm]
&=&  \Pa(x_{t_1}^{a_1})T+ x_{t_1}^{a_1}\Pa(T)\\[0.15cm]
&=& \Pa(x_{t_1}^{a_1})x_{t_2}^{a_2}\cdots x_{t_s}^{a_s}q + x_{t_1}^{a_1} \Pa(x_{t_2}^{a_2}\cdots x_{t_s}^{a_s}q) \\[0.15cm]
&=& \Pa(x_{t_1}^{a_1})x_{t_2}^{a_2}\cdots x_{t_s}^{a_s}q + x_{t_1}^{a_1}\p{\Pa(x_{t_2}^{a_2}\cdots x_{t_s}^{a_s})q+ x_{t_2}^{a_2}\cdots x_{t_s}^{a_s}\Pa(q) } \\[0.15cm]
&=&\p{\Pa(x_{t_1}^{a_1})x_{t_2}^{a_2}\cdots x_{t_s}^{a_s} +  x_{t_1}^{a_1}\Pa(x_{t_2}^{a_2}\cdots x_{t_s}^{a_s})}q+ p\Pa(q) \\[0.15cm]
&=& \Pa(p)q +  p\Pa(q)
\end{eqnarray*}
\end{enumerate}
We will define $\br{\,,\,}_0$ of  recurrently form on the basis of $D$. Put $p\in D$ and define  
\begin{enumerate}
\item $f_p(X)=-g_X(p)$ for $ X\in Mon\{x_1,\cdots,x_n\}$.
\item $f_p(x_1^{-a})=-x_1^{-a}f_{p}(x_1^a)x_1^{-a}$, $a>0$.
\end{enumerate}
And for $x_{l_1}^{b_1}\cdots x_{l_r}^{b_r}\in Mon\{x_1^{\pm},x_2,\cdots,x_n\} $ 
\[
f_p(x_{t_1}^{a_1}\cdots x_{t_n}^{a_n})=f_p(x_{t_1}^{a_1})x_{t_2}^{a_2}\cdots x_{t_n}^{a_n}+ x_{t_1}^{a_1}f_p(x_{t_2}^{a_2}\cdots x_{t_n}^{a_n})
\]
By the universal property of the basis $f$ has a extension to an $A$-homomorphism $f$ in all $D$, now we will take $\br{a,b}_0:=f_a(b)$ and we will prove it is a bracket Poisson on $D$. 
\begin{enumerate}
\item ($\br{a,b}_0+\br{b,a}_0=0 $) Since $\br{a,b}_0$ is a $A$-bilinear function, it is enough to see this when $a,b\in D$ are monomials, put $a=x_1^{-l}a'$ and $b=x_1^{-r}b'$ where $0\leq l,r $,   $a'=x_1^{a_1}\cdots x_n^{a_n}\in B$, and $b'=x_1^{b_1}\cdots x_n^{b_n}\in B$.
\begin{eqnarray*}
\br{a,b}_0+\br{b,a}_0 &:=& f_a(b)+f_b(a)=f_a(x_1^{-r}b')+f_b(x_1^{-l}a') \\[0.15cm]
&:=& \p{f_a(x_1^{-r})b'+x_1^{-r}f_a(b')}+\p{f_b(x_1^{-l})a'+x_1^{-l}f_b(a')}\\[0.15cm]
&:=& \p{-x_1^{-r}f_a(x_1^r)b+x_1^{-r}f_a(b')}+\p{-x_1^{-l}f_b(x_1^l)a+x_1^{-l}f_b(a')}\\[0.15cm]
&:=& \p{x_1^{-r}g_{x_1^{r}}(a)b-x_1^{-r}g_{b'}(a)}+\p{x_1^{-l}g_{x_1^{l}}(b)a-x_1^{-l}g_{a'}(b)}\\[0.15cm]
&=& \p{x_1^{-r}g_{x_1^{r}}(x_1^{-l}a')b-x_1^{-r}g_{b'}(x_1^{-l}a')}+\p{x_1^{-l}g_{x_1^{l}}(x_1^{-r}b')a-x_1^{-l}g_{a'}(x_1^{-r}b')}\\[0.15cm]
&=& x_1^{-r}\p{g_{x_1^r}(x_1^{-l})a'+x_1^{-l}g_{x_1^r}(a')}b-x_1^{-r}\p{g_{b'}(x_1^{-l})a'+x_1^{-l}g_{b'}(a')} \\[0.15cm]
&+& x_1^{-l}\p{g_{x_1^{l}}(x_1^{-r})b'+x_1^{-r}g_{x_1^{l}}(b')}a-x_1^{-l}\p{g_{a'}(x_1^{-r})b'+x_1^{-r}g_{a'}(b')} \\[0.15cm]
&=& x_1^{-r}\p{-x_1^{-l}g_{x_1^r}(x_1^{l})x_1^{-l}a'+x_1^{-l}g_{x_1^r}(a')}b \\[0.15cm]
&-&x_1^{-r}\p{-x_1^{-l}g_{b'}(x_1^{l})x_1^{-l}a'+x_1^{-l}g_{b'}(a')} \\[0.15cm]
&+& x_1^{-l}\p{-x_1^{-r}g_{x_1^{l}}(x_1^{r})x_1^{-r}b'+x_1^{-r}g_{x_1^{l}}(b')}a \\[0.15cm]
&-&x_1^{-l}\p{-x_1^{-r}g_{a'}(x_1^{r})x_1^{-r}b'+x_1^{-r}g_{a'}(b')} \\[0.15cm]
&:=& \p{-x_1^{-l-r}\br{x_1^l,x_1^r}ab+x_1^{-l-r}\br{a',x_1^r}b} \\[0.15cm]
&-&\p{-x_1^{-l-r}\br{x_1^l,b'}a+x_1^{-l-r}\br{a',b'}} \\[0.15cm]
&+& \p{-x_1^{-r-l}\br{x_1^r,x_1^l}ba+x_1^{-r-l}\br{b',x_1^l}a} \\[0.15cm]
&-&\p{-x_1^{-r-l}\br{x_1^r,a'}b+x_1^{-r-l}\br{b',a'}} \\[0.15cm]
&=& \p{-x_1^{-r-l}\br{x_1^l,x_1^r}ab-x_1^{-l-r}\br{x_1^r,x_1^l}ba} \\[0.15cm]
&+&\p{x_1^{-r-l}\br{a',x_1^r}b+x_1^{-r-l}\br{x_1^r,a'}b } \\[0.15cm]
&+& \p{x_1^{-l-r}\br{x_1^l,b'}a+x_1^{-r-l}\br{b',x_1^l}a} \\[0.15cm]
&-& \p{x_1^{-r-l}\br{a',b'}+x_1^{-r-l}\br{b',a'} } \\[0.15cm]
&=&\p{-x_1^{-r-l}\br{x_1^l,x_1^r}ab+x_1^{-l-r}\br{x_1^r,x_1^l}ba} \\[0.15cm]
&=& 0
\end{eqnarray*}
\item  ($\br{a,a}_0=0 $) Let $a=x_1^la'$ where $a'=\sum_t\eta_t x_1^{t_1}\cdots x_n^{t_n}\in B$   and $l\leq 0 $ then 
\begin{eqnarray*}
\br{a,a}_0 &:=&f_a\p{x_1^la'}=\sum_t\eta_tf_a(x_1^l\p{x_1^{t_1}\cdots x_n^{t_n}})\\[0.15cm]
&:=& \sum_t\eta_t\p{f_a(x_1^l)x_1^{t_1}\cdots x_n^{t_n}+x_1^lf_a\p{x_1^{t_1}\cdots x_n^{t_n}} }\\[0.15cm]
&=& f_a(x_1^l)a'+x_1^lf_a(a') \\[0.15cm]
&:=& - x_1^lf_a(x_1^{-l})x_1^la'+x_1^lf_a(a')  \\[0.15cm]
&=& x_1^lg_{x_1^{-l}}(x_1^la')x_1^la'-x_1^lg_{a'}(x_1^la') \\[0.15cm]
&=& x_1 ^l\p{g_{x_1^{-l}}(x_1^l)a'+ x_1^lg_{x_1^{-l}}(a') }x_1^la' - x_1^l\p{ g_{a'}(x_1^l)a'+x_1^lg_{a'}(a')  }\\[0.15cm]
&=& x_1^lg_{x_1^{-l}}(x_1^l)a'a+ x_1^{2l}g_{x_1^{-l}}(a')x_1^la'-  x_1^lg_{a'}(x_1^l)a'-x_1^{2l}g_{a'}(a') \\[0.15cm]
&=&\p{x_1^lg_{x_1^{-l}}(x_1^l)a'a-x_1^{2l}g_{a'}(a')}+\p{x_1^{2l}g_{x_1^{-l}}(a')x_1^la' -  x_1^lg_{a'}(x_1^l)a' } \\[0.15cm]
&=& \p{x_1^{2l}g_{x_1^{-l}}(x_1^{-l})x_1^la'a-x_1^{2l}g_{a'}(a')}+\p{x_1^{2l}g_{x_1^{-l}}(a')x_1^la' +  x_1^{2l}g_{a'}(x_1^{-l})x_1^la' }\\[0.15cm]
&:=&x_1^{2l}\p{\br{x_1^{-l},x_1^{-l}}aa-\br{a',a'}}+x_1^{2l}\p{\br{a',x_1^{-l}}+ \br{x_1^{-l},a'} }a \\[0.15cm]
&=& 0
\end{eqnarray*}

\item ($\br{ab,c}_0=\br{a,c}_0b+a\br{b,c}_0$) It can be proved the same way we proved that $g_X(\cdot)$ is a derivations for all $X\in Mon\{x_1,\cdots,x_n\}$.

\end{enumerate}
\end{proof}

\theorem Let $B=\sigma(A)\langle x_1,\cdots,x_n\rangle$ be a skew PBW extension of $A$ as above and $\br{\,,\,}$ a poisson bracket on $B$, then there exists $\xi\in A$ such that  $\br{a,b}=\xi[a,b]$ for all $a,b\in B$.

\begin{proof}
Since $\br{,}$ is a Poisson bracket  on $B$ then there exists a poisson bracket $ \br{,}_0$ on $\o^{1,n}$ with $\br{,}_0\vert_{B}=\br{,} $ which is A-bilinear function, by theorem (3.5) there exists $\xi\in A$ such that $\br{a,b}_0=\xi[a.b]$ for all $a,b\in \o^{1,n}$ so if $a,b\in B$ then $\br{a,b}=\br{a,b}_0=\xi[a,b]$.

\end{proof}

\section{Some Examples.}
In this section we will show some algebras where we can give a characterization of the poisson brackets.

\begin{enumerate}
\item \textit{The algebra of q-differentiable operators.} $D_{q,h}[x,y]$. Let $q,h\in \K$, $q\neq 0$ consider $\K[y][x;\sigma,\delta]$, $\sigma(y):=qy$ and $\delta(y):=h$. By definitions of skew polynomial ring, it is the $\K$-algebra defined by the relation $xy=qtx+h$. If we put $h=0$ and $q^l- 1\in \K^{\ast}$ for all $l\in\N$ we have that $xr=rx$ and $yr=ry$ for all $r\in \K$.

\item \textit{Algebra of linear partial q-dilatation operators.} For a fixed $q\in \K-\{0\}$, the $\K$-algebra of linear partial q-dilatation operators with polynomial coefficients, respectively, with rational coefficients, is $\K[t_1,\cdots,t_n][H_1^{(q)},\cdots,H_m^{(q)}]$, respectively $\K(t_1,\cdots,t_n)[H_1^{(q)},\cdots,H_m^{(q)}]$ $n\leq m$, subject to the relations:
\begin{eqnarray*}
&t_jt_i=t_it_j \;\;\;\;\;\;\; 1\leq i<j\leq n&\\[0.15cm]
&H_i^{(q)}t_i=qt_iH_i^{(q)} \;\;\;\;\;\;\; 1\leq i\leq n& \\[0.15cm]
&H_j^{(q)}t_i=t_iH_j^{(q)} \;\;\;\;\;\;\; i\neq j& \\[0.15cm]
&H_i^{(q)}H_j^{(q)}=H_j^{(q)}H_i^{(q)} \;\;\;\;\;\;\; 1\leq i\leq n& 
\end{eqnarray*}
If we take $n=m=1$ and $q^l- 1\in\K^{\ast}$ for all $l\in\N$.
\item \textit{Multiplicative analogue of the Weyl algebra.} The $\K$-algebra $O_n(\lambda_{ij})$ is generated by $x_1,\cdots,x_n$ subject to the relations:
\begin{equation*}
x_jx_i=\lambda_{ij}x_ix_j, \;\; 1\leq i<j\leq n
\end{equation*}
where $\lambda_{ij}\in \K-\{0\}$. If we take $n>1$  and $\lambda_{ij}$ as $\N$ independent, i.e. for all $i=1,\cdots,n$ and $m\in \N^n-\{(0,\cdots,0)\ $, $1- \prod_{j\neq i}\lambda_{ij}^{m_j}\in \K^{\ast}$.

\item \textit{3-dimensional skew polynomial algebra \A.} It is given by the relations 
\begin{equation*}
yz-\alpha zy=\lambda, \;\; zx-\beta xz=\mu,\;\; xy-\gamma yx=v
\end{equation*}
such that $\lambda,\mu,v\in \K+\K x+\K y+\K z$, and $\alpha,\beta, \gamma\in \K-\{0\}$. If we take $\lambda,\mu =0 $ and $\alpha,\beta,\gamma$, $\N$-independent.

\item \textit{Quantum Space $S_q$.} Let $\K$ be a commutative ring and let $\textbf{q}=[q_{ij}]$ be a matrix with entries in $\K^{\ast}$ such that $q_{ii}=1=q_{ij}q_{ji}$ for all $1\leq i,j \leq n$. The $\K$-algebra $S_q$ is generated by $x_1,\cdots,x_n$, subject to the relations
\begin{equation*}
x_ix_j=q_{ij}x_jx_i
\end{equation*}
If we take $n>1$ and $q_{ij}$ $\N$-independent.

\item \textit{Witten's deformation of $\mathcal{U}(\mathcal{SL}(2,\K))$.} Let $\underline{\xi}=(\xi_1,\cdots,\xi_7)$ a 7-tuple of parameters, it is generated by $x,y,z$ subject to relations 
\begin{equation*}
xz-\xi_1 zx=\xi_2x, \;\; zy-\xi_3yz=\xi_4y,\;\; yx-\xi_5xy=\xi_6z^2+\xi_7z
\end{equation*}
If we take $\xi_7,\xi_6,\xi_2=0$ and $\xi_1,\xi_3,\xi_5\in \K^{\ast} $ and $\N$-independent.
\item \textit{Quantum symplectic space.} $O_q(\mathcal{SP}(\K^{2n}))$. For every nonzero element $q\in \K$, one defines this quantum algebra $O_q(\mathcal{SP}(\K^{2n}))$ to be the algebra generated by $\K$ and the variables $y_1\cdots,y_n,x_1,\cdots,x_n$ subject to the relations 
\begin{eqnarray*}
&y_jx_i=q^{-1}x_iy_j, \;\;\; y_iy_j=y_iy_j,\;\;\; 1\leq i<j\leq n&\\[0.15cm]
&x_jx_i=q^{-1}x_ixj, \;\;\; x_jy_i=qy_ix_j,\;\;\; 1\leq i<j\leq n&\\[0.15cm]
&x_iy_i-q^2y_ix_i=(q^2-1)\sum_{l=1}^{i-1}q^{i-l}y_lx_l, \;\;\; 1\leq i\leq n&\\[0.15cm]
\end{eqnarray*}
If we take $n=1$ and $q$ $\N$-independent.
\end{enumerate}
\section*{ Acknowledgements}
The author would like to thank 
Jose Oswaldo Lezama Serrano for many helpful suggestions which were so important to develop the article.

\bibliographystyle{alpha}

\end{document}